\documentclass[aip, cha, reprint]{revtex4-1}
\usepackage{amsmath,amssymb,amsthm}
\newcommand{\Mod}[1]{\ (\text{mod}\ #1)}
\usepackage{graphicx}% Include figure files
\usepackage{dcolumn}% Align table columns on decimal%point
\usepackage{bm}% bold math
\usepackage{subfig}% subfigs
\usepackage{enumerate}
\usepackage{color,transparent} % color and transparent
\usepackage{bbm}
\usepackage{multirow}

\graphicspath{{./newfigures/}}
\begin{document}

\preprint{AIP/123-QED}

\title{A rough-and-ready cluster-based approach for extracting finite-time coherent sets from sparse and incomplete trajectory data}

\author{Gary Froyland}
\affiliation{
School of Mathematics and Statistics, The University of
New South Wales, Sydney NSW 2052, Australia
}%
\author{Kathrin Padberg-Gehle}
\affiliation{
Technische Universit\"{a}t Dresden,
Fachrichtung Mathematik,
Institut f\"{u}r Wissenschaftliches Rechnen,
D-01062 Dresden, Germany
}

\date{\today}

\begin{abstract}
We present a numerical method to identify regions of phase space that are approximately retained in a mobile compact neighbourhood over a finite time duration.
Our approach is based on spatio-temporal clustering of trajectory data.
The main advantages of the approach are the ability to produce useful results (i) when there are relatively few trajectories and (ii) when there are gaps in observation of the trajectories as can occur with real data.
The method is easy to implement, works in any dimension, and is fast to run.
\end{abstract}

\pacs{}
\keywords{coherent sets, Lagrangian coherent structure, spatio-temporal clustering.}
\maketitle

\begin{quotation}

Coherent features in time-dependent dynamical systems are difficult to identify, and considerable effort has been put into the development of identification algorithms.
Most approaches require knowledge of the dynamical system or high-resolution trajectory information, which in applications may not be available.
We present a trajectory-based method that is aimed squarely at the situation where the available information is poor:  there are few trajectories, the available trajectories do not span the full time duration under consideration, and there are missing observations within trajectories.
As our method is very simple to implement and fast to run, it also provides a rapid ``first cut'' coherent structure analysis even in situations where the full dynamical system or high-resolution trajectory data is available.
 \end{quotation}

 \section{Introduction}
There are a number of different concepts that describe the notion of coherent behaviour in time-dependent dynamical systems.
Probabilistic approaches define finite-time coherent sets\cite{FSM10,F13,FPG14} as regions of phase space that minimally mix with the surrounding phase space during a specified time duration of finite length.
Lagrangian coherent structures can be defined as material lines that extremize a certain stretching or shearing quantity\cite{haller_11,haller_12}, while another approach tries to identify curves on which local dynamics approximates local rigid-body motion\cite{bollt2014}.
There are also topological\cite{AT12} and ergodicity-based\cite{BM12} descriptions of coherence, although these are not designed for aperiodic dynamics.
Finally, a recent geometric characterisation\cite{F15} defines finite-time coherent sets as those sets with boundary to volume ratios that remain minimal under the evolution of the dynamics, and proves that such a characterisation arises naturally as the advective limit of the probabilistic approaches\cite{FSM10,F13}.

In the present paper we develop cluster-based techniques to highlight distinct groups of trajectories that remain in compact, approximately spherical subregions of phase space over a finite time duration.
Let $\Phi:\mathbb{R}^d\times [0,T]\to \mathbb{R}^d$ denote the flow of a continuous time dynamical system on $\mathbb{R}^d$, i.e. $\Phi(x,t)$ denotes the state of the system at time $t$ with initial value $x$ (at time 0).
We define a dynamic metric
\begin{equation}
\label{dynmetric}
\mathbf{D}(x,y):=\int_0^T\rho(\Phi(x,t),\Phi(y,t))^2 dt,
\end{equation}
based on some metric $\rho$ on $\mathbb{R}^d$.
For example, if $\rho$ is the Euclidean metric, then $x$ and $y$ are close according to $\mathbf{D}$ provided they remain close in a Euclidean sense averaged over the time interval $[0, T]$.

In this general setup, one is free to choose $\rho$ and also how the terms $\rho(\Phi(x,t),\Phi(y,t))$, $t\in [0, T]$ are combined to form $\mathbf{D}(x,y)$.
Setting $\rho$ to be the Euclidean metric is a natural choice if the trajectory data lies in $\mathbb{R}^d$ and shortly we will give geometric reasons for why this is a good choice.
The sum-of-squares combination is a convenient form for the specific numerical clustering approach proposed below.
One could alternatively define, for example, $\mathbf{D}(x,y):= \int_0^T\rho(\Phi(x,t),\Phi(y,t))^p dt$, for $1\le p<\infty$, or $\mathbf{D}(x,y):=\max_{t\in[0, T]} \rho(\Phi(x,t),\Phi(y,t))$.

In practice, suppose we have $n$ trajectories given at discrete times $\{0, 1,\ldots, T\}$, denoted $x_{i,t}\in \mathbb{R}^d, i=1,\ldots,n, t=0,\ldots,T$.
We wish to cluster the initial points $x_{i,0}\in\mathbb{R}^d$ according to $\mathbf{D}$.
The discrete-time version of (\ref{dynmetric}) is
\begin{equation}
\label{dynmetric2}
\mathbf{D}(x_{i,0},x_{j,0})=\sum_{t=0}^{T} \rho(x_{i,t},x_{j,t})^2,
\end{equation}
for $1\le i,j\le n$.
At this point, one could calculate $n(n-1)/2$ interpoint distances $\mathbf{D}(x_{i,0},x_{j,0})$, $1\le i<j\le n$.
This general approach of clustering using the dynamic metric (\ref{dynmetric}) or (\ref{dynmetric2}) is very flexible and in principle one could employ any suitable (according to the three properties outlined below) clustering method on $\mathbb{R}^d$ from the vast number available (see e.g.\ Ref.~\onlinecite{rokach_2010}):  centroid-based algorithms like k-means\cite{lloyd1982,arthur_vassil} and fuzzy c-means\cite{bezdekbook,bezdek87}, spectral methods
 \cite{fiedler1,shi}, density-based clustering \cite{Ester_et_al1996,ankerst_1999}, and methods based on community detection \cite{fortunato} (for example, modularity
  \cite{newman}). 
  %A recent survey of clustering algorithms can be found in Ref.~\onlinecite{rokach_2010}.

The sum-of-squares form of $\mathbf{D}$ and choice of Euclidean metric for $\rho$ allows us to further rewrite (\ref{dynmetric2}) as \begin{equation}
\label{dynmetric3}
\mathbf{D}(x_{i,0},x_{j,0})=\sum_{t=0}^{T} \rho(x_{i,t},x_{j,t})^2=\|X_i-X_j\|^2,
\end{equation}
for $1\le i,j\le n$, where $X_i=(x_{i,0},x_{i,1},\ldots,x_{i,T})$.
Thus, we have a convenient representation of the dynamic metric $\mathbf{D}$ on $\mathbb{R}^d$ as the squared Euclidean distance between trajectories in $\prod_{t=0}^{T} \mathbb{R}^d=\mathbb{R}^{d(T+1)}$.

As elaborated in the next section we use the fuzzy $c$-means clustering algorithm\cite{bezdekbook,bezdek87} on $\mathbb{R}^{d(T+1)}$.
Our reasons for using fuzzy $c$-means are threefold (but not necessarily exclusive to fuzzy $c$-means): \begin{enumerate}
 \item When searching for $K$ clusters\footnote{If $K$ is not known \emph{a priori}, to determine $K$ one can employ any of several existing methods which automatically produce an optimal number of clusters from datasets.}, fuzzy $c$-means produces $K$ auxilliary ``centres'' and aims to allocate data to clusters by reducing the total squared distance from the data to their corresponding centre.  If $\rho$ is the Euclidean metric, then fuzzy $c$-means will favour clusters that are close to spherical at each time instant.
     Such clusters therefore will not ``spread out'' in phase space, will remain in an approximately tubular region in lower-dimensional space-time (phase space plus one time coordinate, see Figure \ref{fig:fc_visualize}), and will on average have low boundary size to volume at each time instant.
      The low boundary size relative to volume property is compatible with probabilistic approaches \cite{FSM10,F13} and geometric approaches \cite{F15} to finite-time coherent sets.
   \item  Fuzzy $c$-means provides feedback in the form of the membership value describing the likelihood that a trajectory belongs to a cluster.  Because finite-time coherent sets do not necessarily fully partition the phase space, we can identify non-coherent collections of trajectories as those with a low membership for all clusters.
   \item Fuzzy $c$-means is computationally efficient, particularly for large numbers of trajectories.
\end{enumerate}
Finally, we note that once the initial points $\{x_{i,0}\}_{1\le i\le n}$ have been clustered, the full trajectories are also clustered, as by definition trajectories remain within the same cluster for all $t=0,\ldots,T$.

Clustering trajectory data is a recent problem in the analysis of spatio-temporal datasets, with many contributions found in the data mining literature. We refer the reader to Ref.~\onlinecite{KMNR10} for a recent review and to the literature review in Ref.~\onlinecite{izakian_13} for a brief summary of the different approaches used for spatio-temporal clustering. Ref.~\onlinecite{izakian_13} proposes an augmented fuzzy c-means algorithm, with different weights for the temporal and spatial components.
Ref.~\onlinecite{nanni_2006} and the thesis \onlinecite{nanni_2002} introduce a distance measure essentially identical in form to (\ref{dynmetric}).  Ref.~\onlinecite{nanni_2006} uses this metric with density-based algorithms\cite{ankerst_1999}  to cluster trajectories in geo-referenced data sets and to find optimal time intervals for clustering.
The papers \onlinecite{nanni_2006,izakian_13,nanni_2002} do not consider how to treat incomplete trajectory data. Other distance measures have been proposed to account for application-specific purposes, e.g. for studying movement patterns in traffic\cite{pelekis_2007,Rinzivillo_2008}. However, to the best of our knowledge, spatio-temporal clustering approaches have not been employed for studying transport phenomena and coherent behaviour in time-dependent dynamical systems and an exploration of the tuning of clustering methods to this application has not been undertaken.

Our main contributions are (i) posing the problem of identifying finite-time coherent sets as an objective trajectory-based clustering problem, (ii) developing a methodology for handling incomplete data that consistently uses all available data, and (iii) indicating some rules of thumb for applying these techniques in practice.

An outline of the paper is as follows.
We describe our approach first in the situation where there are a finite number of trajectories available, sampled at a finite number of times.
We then consider combinations of a continuum of trajectories and a continuum of observation times.
Section \ref{sec:fulldata} concludes by showing that our clustering framework for coherent sets is objective and independent of isotropic scaling of space and time.
Our method handles missing trajectory data naturally and we discuss this in Section \ref{sec:missingdata}.
A discussion of false positives, possibly inaccurate results, and how to identify these is in Section \ref{sec:gowrong};  we also outline some rules of thumb for parameter choices.
Section \ref{sec:examples} illustrates the approach for several examples:  firstly in one-dimensional dynamics, where the geometry of the spatio-temporal clustering is more transparent, secondly for the well-known double-gyre flow and the transitory double gyre flow for comparison with existing coherent set identification approaches, and thirdly on ocean surface drifter data.
We demonstrate how reasonable results can be achieved even when large percentages of trajectory observations are missing, and when the trajectory dataset is comprised of trajectories much shorter than the full time duration under analysis.

\section{Full data case}\label{sec:fulldata}

We first describe the case where all trajectories span the finite time duration and there are no missing observations.
We begin by describing our setup in the situation where there are a finite number of trajectories sampled at a finite collection of time instances, this means that all trajectories are sampled at all time instances;  we then follow with versions that are continuous in space and/or time.

\subsection{Discrete setting: a finite number of finitely-sampled trajectories.}
\label{sect:discrete}

Suppose we have $n$ trajectories of maximum length $T+1$, denoted $x_{i,t}\in \mathbb{R}^d, i=1,\ldots,n, t=0,1,\ldots,T$.
We consider each trajectory $\{x_{i,t}\}_{0\le t\le T}$ as a point $X_i=(x_{i,0},x_{i,1},\ldots,x_{i,T})\in \mathbb{R}^{d(T+1)}$, which we also refer to as a trajectory.
We imagine $\mathbb{R}^{d(T+1)}$ as $\prod_{t=0}^T \mathbb{R}^d$, where the product of copies of our phase space $\mathbb{R}^d$ is ordered in increasing time $t$.
The fuzzy $c$-means clustering algorithm\cite{bezdekbook,bezdek87} is a soft clustering based on the calculation of a centre for each cluster and a likelihood of membership of each data point to each centre.
Suppose we have identified $K$ cluster centres $C_k\in\mathbb{R}^{d(T+1)}, k=1,\ldots,K$.
We may decompose the $C_k$ as $(c_{k,0},c_{k,1},\ldots,c_{k,T})\in \mathbb{R}^{d(T+1)}$, so that each $c_{k,t}\in \mathbb{R}^d$ can be regarded as a point in phase space $\mathbb{R}^d$ at time $t$.
Associated with each trajectory $X_i$ is a likelihood $0\le u_{k,i}\le 1$ of $X_i$ being associated with the cluster centre $C_k$.

Given trajectories $X_i, i=1,\ldots,n$, cluster centres $C_k$, $k=1,\ldots,K$, and membership likelihoods $u_{k,i}, i=1,\ldots,n, k=1,\ldots,K$, the total ``goodness of fit'' of the memberships of trajectories in clusters is measured by likelihood-weighted intracluster distances, which we wish to minimise:
\begin{equation}
\label{obj1}
\sum_{k=1}^K\sum_{i=1}^n u_{k,i}^m\|X_i-C_k\|^2=\sum_{k=1}^K\sum_{i=1}^n u_{k,i}^m\sum_{t=0}^T\|x_{i,t}-c_{k,t}\|^2.
\end{equation}
This minimisation is subject to the constraints that (i) $\sum_{k=1}^K\ u_{k,i}=1$ for $i=1, \ldots, n$ and (ii)
$u_{k,i}\ge 0$ for all $k=1,\ldots,K$, $i=1,\ldots,n$.
The parameter $m>1$ is the fuzziness exponent.
Increasing $m$ corresponds to softer clusters, while as $m$ approaches 1, the membership likelihoods converge to either $0$ or $1$, resulting in a hard clustering\cite{bezdek87} (this latter effect is most easily seen from the update rule (\ref{uupdateinc}) below).
The basic fuzzy $c$-means algorithm\cite{bezdekbook,bezdek87} in our notation above proceeds as follows.

\vspace{.2cm}
\noindent\textbf{Algorithm 1:}
\begin{enumerate}
\item Initialize membership values $u_{k,i}$ either randomly or computed via step 3 based on an initial seeding of $K$ centres (e.g. randomly or by the $k$-means++ algorithm \cite{arthur_vassil})
\item Calculate centres:
\begin{equation}
\label{cupdateinc}
C_k=\frac{\sum_{i=1}^n u_{k,i}^mX_i}{\sum_{i=1}^n u_{k,i}^m},
\end{equation}
$k=1,\ldots,K$.
\item Update membership values:
\begin{equation}
\label{uupdateinc}
u_{k,i}=\frac{1/\|{X}_i-C_k\|^{2/(m-1)}}{\sum_{j=1}^K \left(1/\|{X}_i-C_j\|^{2/(m-1)}\right)},
\end{equation}
$k=1,\ldots,K$, $i=1,\ldots,n$.
\item Evaluate objective (\ref{obj1}).  If the improvement in the objective is below a threshold, go to step 5;  otherwise go to step 2.
\item Output cluster centres $C_k\in \mathbb{R}^{d(T+1)}, k=1,\ldots,K$ and membership likelihoods $u_{k,i}\in [0,1], k=1,\ldots,K$, $i=1,\ldots,n$.
\end{enumerate}
The update rule in step 2 is constructed by fixing the $u_{k,i}$ and choosing $C_k$ so that the gradient of the objective (\ref{obj1}) is zero.
Similarly, the update rule in  step 3 is constructed by fixing the $C_k$ and choosing the $u_{k,i}$  so that the gradient of the Lagrangian incorporating (\ref{obj1}) and the constraint $\sum_{k=1}^K u_{k,i}=1$ is zero.
\vspace{.2cm}

\noindent\textbf{Implementation in MATLAB:} If $X$ is an $n\times d(T+1)$ array of $n$ trajectories in $\mathbb{R}^d$ (i.e. the rows of $X$ are the vectors $X_i$, $i=1,\ldots,n$ discussed above),  Algorithm 1 is implemented in MATLAB by the function {\tt fcm} in the Fuzzy Logic Toolbox:
\begin{verbatim}
opts(1)=m;
[c,u]=fcm(X,K,opts);
\end{verbatim}
If {\tt fcm} is called without {\tt opts}, the default value of $m$ is 2.
For $d=2$, to display the membership values for cluster $k$ at time slice $t\in\{0,1,\ldots,T\}$, one can use
\begin{verbatim}
scatter(X(:,2*t+1),X(:,2*t+2),[],u(k,:),'.');
\end{verbatim}
For $d=3$, one can similarly use ${\tt scatter3}$.

\begin{figure}[hbt]
  \centering
  \includegraphics[width=8cm, clip=true]{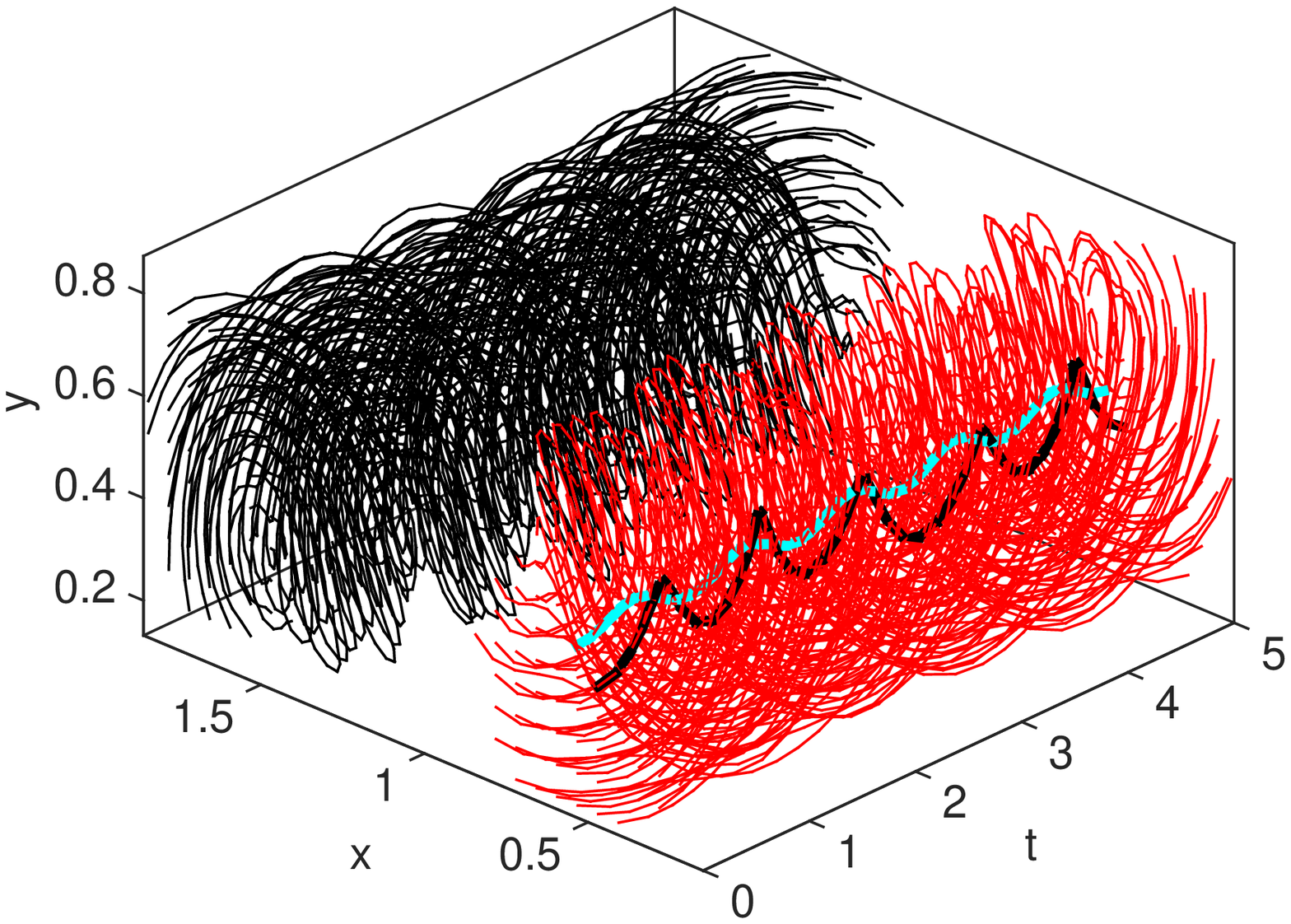}\\
  \caption{Result of clustering of the double-gyre flow with $n=512$, $K=2$, $m=2$, $d=2$, $T=50$ (corresponding to the real-time interval $[0, 5]$ in time steps of $0.1$).  Trajectories $X_i$ with membership values $u_{k,i}>0.9$, $k=1,2$ are shown.  Note that in this figure we display space-time using time as a third coordinate, whereas the clustering computations in Algorithm 1 occur in $\mathbb{R}^{102}$.  One sees that the red and black trajectories remain in a compact region as time evolves. The probabilistic centre of the cluster corresponding to the red trajectories is shown as a meandering black curve, the light blue curve is the centre as computed by Algorithm 1.}\label{fig:fc_visualize}
\end{figure}

We remark that the centres $C_k=(c_{k,0},c_{k,1},\ldots,c_{k,T})$, $k=1,\ldots,K$, are generally not true trajectories of the dynamical system, although they may be remarkably close to true trajectories in some cases.
For each $k=1,\ldots,K$, one can identify the maximum likelihood trajectory $X_{i_k^*}$ for the $k^{\rm th}$ cluster, where $i_k^*=\arg \max_i u_{k,i}$.
The trajectory $X_{i^*_k}$ is the most likely to belong to the $k^{\rm th}$ cluster and may be thought of as a ``probabilistic centre'' of the cluster.
The probabilistic centers can be interpreted as ``low dimensional representations'' of the macroscopic behavior of the system, as they describe the coherent motion of trajectories in the corresponding cluster.
We illustrate both the centre and maximum likelihood trajectory in Figure \ref{fig:fc_visualize}.

Note that there is potential for one to include weights as coefficients for the terms $\|x_{i,t}-c_{k,t}\|^2$ in (\ref{obj1}), which could depend on $i,t,$ or $k$.
A particularly important example is the inclusion of weights $q_i\ge 0$ corresponding to the ``mass'' assigned to a point $x_{i,0}$.
For example, if one is searching for coherent regions in an oil or chemical spill in the ocean, one is likely interested in the behavior of the oil or chemical, rather than the water.
In order to obtain clusters that focus on the nonuniform distribution of oil or chemical, one can replace $u_{k,i}^m$ with $q_iu_{k,i}^m$ in (\ref{obj1}) and (\ref{cupdate}).

Also note that at present, clustering into spheres is preferred by the Euclidean norm.
If one wishes to favour clustering into ellipsoids, with orthogonal semi-axis vectors $v_1,\ldots,v_d$ and corresponding semi-axis lengths $\ell_1,\ldots,\ell_d$, then one may simply scale the $x_{i,t}$ data in $\mathbb{R}^d$ along each $v_j$ by $1/\ell_j$, $j=1,\ldots,d$, use the Euclidean norm in the objective above, and then rescale the data along each $v_j$ by $\ell_j$, $j=1,\ldots,d$.
Other distance functions could be used to replace Euclidean distance, but in the absence of specific replacement motivations based on known properties of the underlying dynamical system, Euclidean distance represents a natural isotropic default distance metric.

Sections \ref{sect:semicts1}--\ref{sect:fullycts} outline extensions of the above setup to situations where either one or both of the spatial data or temporal data are on a continuum.
These constructions are mainly of a theoretical nature, but have been included to (i) demonstrate what the analogous objects are in a continuum setting if e.g.\ a full dynamical systems model were available, and (ii) indicate how the discrete ``finite data'' setting above is a special case (constructed by subsampling in space and/or time) of the continuum ``full model'' setting.
Sections \ref{sect:semicts1}--\ref{sect:fullycts} could be omitted on a first reading.

\subsection{Semi-continuous setting \#1: a continuum of initial points, with trajectories finitely-sampled in time.}
\label{sect:semicts1}

Suppose we have a continuum of initial points in a set
$A\subset \mathbb{R}^d$.
We now write trajectories as $x(x_0,t)\in \mathbb{R}^d, x_0\in A, t=0,1,\ldots,T$.
Individual trajectories $\{x(x_0,t)\}_{0\le t\le T}$ for fixed $x_0\in A$ are still regarded as elements of $\mathbb{R}^{d(T+1)}$ as before, and we write an individual trajectory as $X(x_0)\in \mathbb{R}^{d(T+1)}$.
Note that  the likelihoods $u_{x_0,k}$ are also continuously parameterised by $x_0\in A$, and we write these in functional form as $u_k(x_0)$, so that $u_k:A\to [0,1]$, $k=1,\ldots,K$.
In the finite-trajectory setting, the initial points of trajectories need not be uniformly distributed over the phase space, nor be given a uniform weight.
If one  wishes to model the evolution of a passive tracer field with nonuniform density, one will either have a greater density of points in areas of high tracer density or apply weights to points with higher tracer concentration.
To capture this effect in the continuum setting, we need a density function $q:A\to \mathbb{R}$, satisfying $\int_{A}q(y)\ dy=1$.
We interpret $\int_{B_\epsilon(x_0)}q(y)\ dy$ as the fraction of initial points that belong to an $\epsilon$-neighbourhood of $x_0$.
For example, if the initial $x_0$ are uniformly sampled over $A$, then $q\equiv 1/{\rm vol}(A)$.

Equation (\ref{obj1}) now reads
\begin{eqnarray}
\nonumber\lefteqn{\sum_{k=1}^K \int_A u_{k}(x_0)^m\|X(x_0)-C_k\|^2q(x_0)\ dx_0}\\
\label{obj2}&=&\sum_{k=1}^K \sum_{t=0}^T \int_A  u_k(x_0)^m\|x(x_0,t)-c_{k,t}\|^2q(x_0)\ dx_0
\end{eqnarray}

\begin{figure}[hbt]
\centering
\begin{tabular}{cc}
  \includegraphics[width=0.45\columnwidth, clip=true]{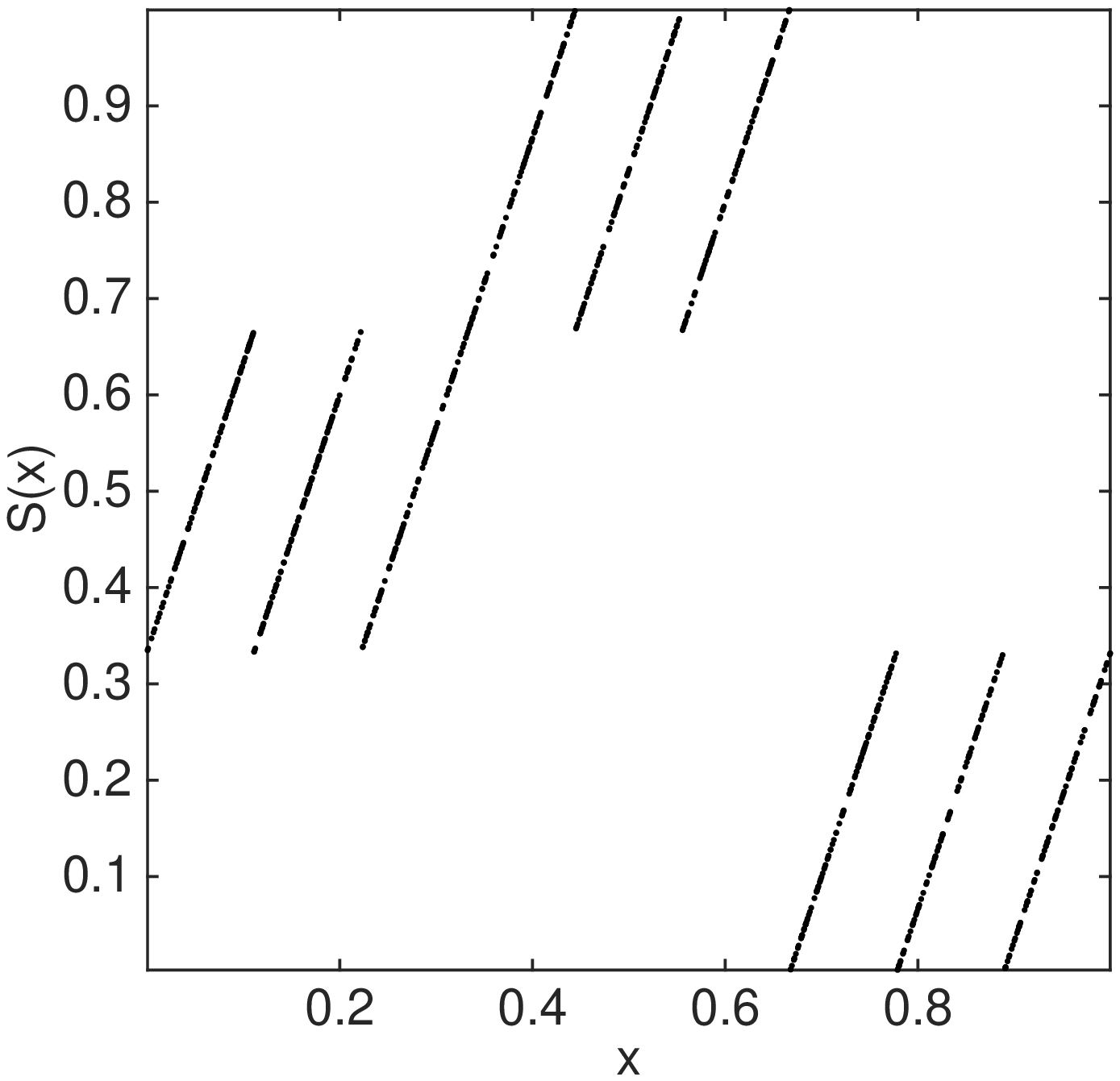}&
\includegraphics[width=0.45\columnwidth, clip=true]{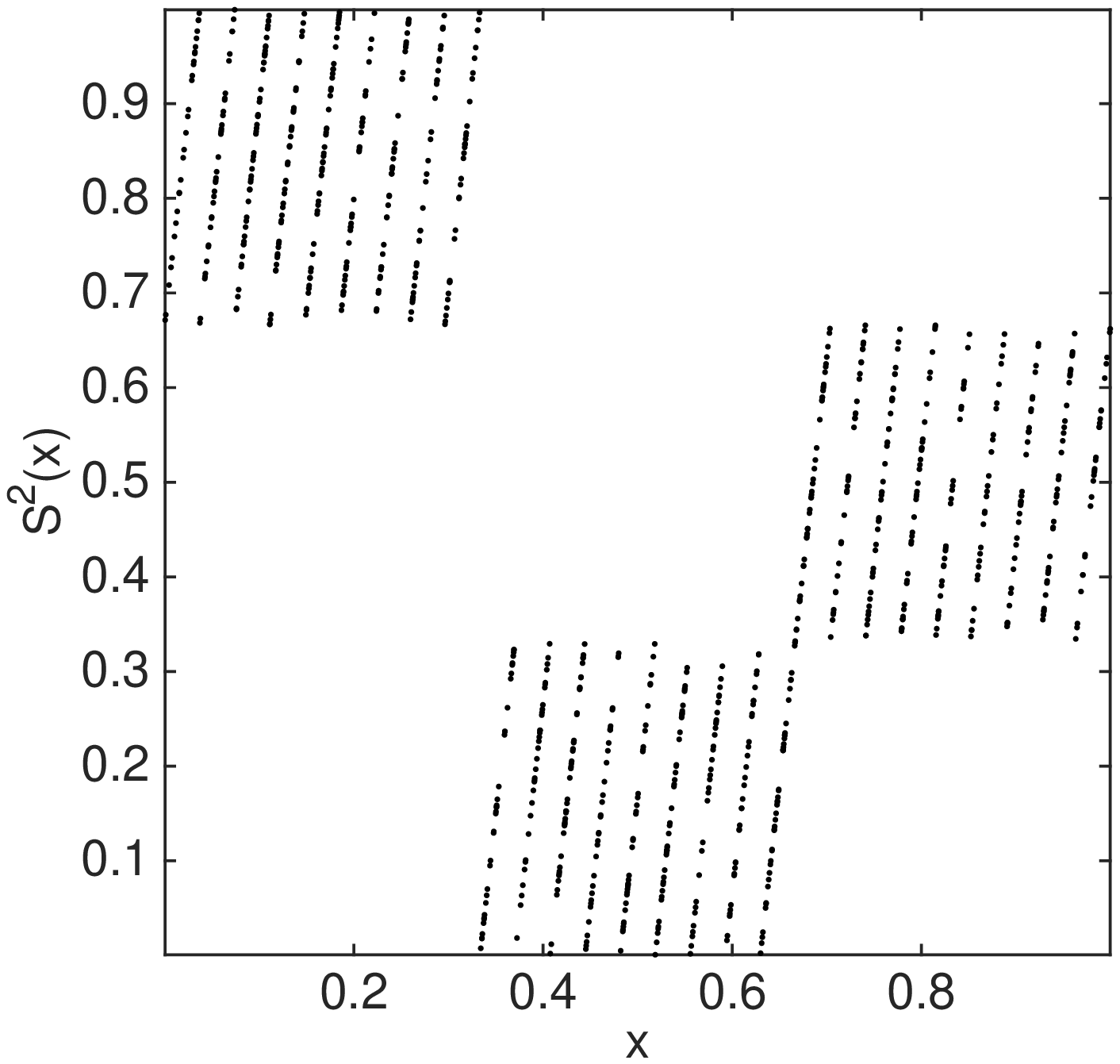}\\
 {\scriptsize (a)} & {\scriptsize (b)}\\
 \end{tabular}
  \caption{Graphs of the interval map $S$ (\ref{eq:3map}), which permutes the intervals $[0,1/3],[1/3,2/3],[2/3,1]$. (a) $x$ vs $S(x)$; (b) $x$ vs $S^2(x)$. Considered as a time series of length 2, Algorithm 1 would seek clusters in these graphs, considered as a subset of $[0,1]^2$.}\label{fig:intervalmap}
\end{figure}

Here is a simple example to help visualise what is going on.
Let phase space be $[0,1]$, and consider trajectories of length two, generated by a map $S:[0,1]\circlearrowleft$.
Geometrically, we look for clusters in data of the form $(x(x_0,0),x(x_0,1))$, for all $x_0\in [0,1]$, which is nothing but the (weighted, if $q$ is not constant) graph of $S$ considered as a one-dimensional subset of $[0,1]^2$;  see Figure \ref{fig:intervalmap}(a).  
Figure \ref{fig:intervalmap}(b) shows clusters in data of the form $(x(x_0,0),x(x_0,2))$.

\subsection{Semi-continuous setting \#2: a finite number of continuously-sampled trajectories.}
\label{sect:semicts2}

Suppose now that instead of a continuum of finitely-sampled trajectories, we now have a finite collection of trajectories observed continuously in time: $x_i(t), i=1,\ldots,n, t\in [0,T]$, with optional weights $q_i$, $i=1,\ldots,n$, representing the mass of the point $x_i(0)$.
We consider $x_i:[0,T]\to \mathbb{R}^d$, $i=1,\ldots,n$ as a finite number of continuous mappings from $[0,T]$ to $\mathbb{R}^d$.
Similarly, the cluster centres $c_k:[0,T]\to \mathbb{R}^d$, $k=1,\ldots,K$ are a finite number of (not necessarily continuous) mappings from $[0,T]$ to $\mathbb{R}^d$.
The RHS of (\ref{obj1}) becomes
\begin{equation}
\label{obj3}
\sum_{k=1}^K \sum_{i=1}^n u_{k,i}^m\int_0^T\|x_i(t)-c_k(t)\|^2q_i\ dt.
\end{equation}

To visualise the integral, imagine we have a one-dimensional (time-dependent) flow in a phase space $\mathbb{R}$, and $x_i(t)$ is one trajectory from $t=0$ to $t=T$.
The integral (\ref{obj3}) computes the total ($q$-weighted) squared distance between the graphs of the functions $x_i$ and $c_k$ in $[0,T]\times \mathbb{R}$;  or in other words, the squared $L^2$ distance $\|x_i-c_k\|^2_{L^2([0,T])}$.

\subsection{Fully continuous setting: a continuum of initial points, continuously-sampled trajectories.}
\label{sect:fullycts}
Combining the constructions from the previous two paragraphs, we now have functions $x:A\times [0,T]\to \mathbb{R}^d$.
Our likelihood functions remain as $u_k:A\to \mathbb{R}^+$, $k=1,\ldots,K$.
The RHS of (\ref{obj1}) becomes
\begin{equation}
\label{obj4}
\sum_{k=1}^K \int_0^T\left(\int_A u_k(x_0)^m\|x(x_0,t)-c_k(t)\|^2q(x_0)\ dx_0\right)\ dt.
\end{equation}

\subsection{Isotropic scaling of space and time has no effect}

Given that our clustering is occurring in the product space formed from as many copies of our phase space as there are time instants, it is pertinent to consider the effect, if any, of isotropically scaling space and time.
We show that in fact, there is no real effect caused by such scaling.

In the fully discrete setting, by (\ref{obj1}) if space were scaled isotropically by a factor $\alpha$ and time by a factor $\beta$, then (\ref{obj1}) would simply increase by $\alpha^2$.
Thus the cluster centres and likelihood functions that minimise (\ref{obj1}) are simply isotropically scaled and unchanged, respectively, under this isotropic scaling of space and/or time.

In the continuum setting, we again consider scaling space isotropically by a factor $\alpha$ and time by a factor $\beta$.
This amounts to defining new ``primed'' variables: $x'_0=\alpha x_0$, $t'=\beta t$, $A'=\alpha A$, $T'=\beta T$,  $x'(x'_0,t')=\alpha x(x_0,t), c'_k(t')=\alpha c_k(t), u'_k(x'_0)=u_k(x_0)$, and $q'(x'_0)=q(x_0)/\alpha^d$.
Then changing variables from $x_0$ to $x'_0$ and from $t$ to $t'$ we have
\begin{eqnarray*}
\lefteqn{(\ref{obj4})}\\
&=&\sum_{k=1}^K \int_0^{\beta T}\left(\int_{\alpha A} u_k(x'_0)^m\|x(x'_0,t')-c_k(t')\|^2\right.\\
&&\left.\qquad\qquad\qquad q(x'_0)\ \alpha^{-d}\beta^{-1}\ dx'_0\right)\ dt'\\
&=&\alpha^{-2}\beta^{-1}\sum_{k=1}^K \int_0^{T'}\left(\int_{A'} u'_k(x'_0)^m\|x'(x'_0,t')-c'_k(t')\|^2\right.\\
&&\qquad\qquad\qquad\left. q'(x'_0)\ dx'_0\right)\ dt'.
\end{eqnarray*}
Thus, switching to the primed coordinates will simply increase the objective (\ref{obj4}) by a constant factor $\alpha^2\beta$ over the original unprimed coordinates.
Again, the cluster centres and likelihood functions that minimise (\ref{obj4}) are isotropically scaled and unchanged, respectively, under this isotropic scaling of space and/or time.
In particular, the clustering algorithm does not care how space is scaled against time.

\subsection{Frame-independence}

To check frame-independence of an algorithm, one applies the algorithm to an original dataset, then subjects the dataset to a (possibly time-dependent) affine transformation, where the linear part is orthogonal.
If the algorithm applied to the transformed dataset yields the transformed output of the original dataset, then the algorithm is frame-independent; see Ref.~\onlinecite{truesdellnoll} for details.

We consider the situation where we have a finite collection of finitely-sampled trajectories;  the arguments presented apply equally to the other situations discussed in Sections \ref{sect:semicts1}--\ref{sect:fullycts}.

\vspace{.2cm}
\noindent\textbf{Proposition:}
Algorithm 1 is frame-independent.
\vspace{.2cm}

\noindent\emph{Proof:}
Let $\{x_{i,t}\}_{1\le i\le n,0\le t\le T}$ be an original collection of trajectories.
Apply Algorithm 1 to $\{x_{i,t}\}$ to obtain centres $C_k=(c_{k,0},\ldots,c_{k,T})\in \mathbb{R}^{d(T+1)}$, $k=1,\ldots,K$ and likelihoods $u_{k,i}, k=1,\ldots,K,$ $i=1,\ldots,n$ that minimise (\ref{obj1}).
Denote the transformed trajectories $y_{i,t}:=O_{t}x_{i,t}+o_t$, where $O_t$ is an orthogonal $d\times d$ matrix and $o_t\in \mathbb{R}^d$.
Form transformed centres ${c}'_{k,t}:=O_tc_{k,t}+o_t$.
Notice that (\ref{obj1}) has the same value when evaluated with $\{x_{i,t}\}$, $\{c_{k,t}\}$ and $\{u_{k,i}\}$, and with $\{y_{i,t}\}$, $\{c'_{k,t}\}$ and $\{u_{k,i}\}$.
This is because the transformation $x\mapsto O_tx+o_t$ is an isometry with respect to the Euclidean norm for each $t=0,1,\ldots,T$.
Because $\{c_{k,t}\}$ and $\{u_{k,i}\}$ minimise (\ref{obj1}) for the dataset $\{x_{i,t}\}$, one has $\{c'_{k,t}\}$ and $\{u_{k,i}\}$ minimise (\ref{obj1}) for the dataset $\{y_{i,t}\}$.\ $\square$

If we use a non-standard inner product $\langle\cdot,\cdot\rangle':=x^\top Qx$ for some symmetric positive-definite $d\times d$ matrix $Q$ to define a norm $\|\cdot\|'$ on each phase space slice $\mathbb{R}^{d}$, then an analogous proposition would hold under transformations $x\mapsto O_tx+o_t$ provided $O_t^*=Q^{-1}O_t^\top Q=O_t^{-1}$.

\section{Treatment of missing data}\label{sec:missingdata}

Missing data can be treated naturally in our spatio-temporal clustering framework.
Taking the finitely sampled, finite trajectory setting of Section \ref{sect:discrete}, by \emph{missing data}, we mean a trajectory $\{x_{i,t}\}_{0\le t\le T}$ where the values $x_{i,t}$ are available only on a strict subset of time instances $\mathcal{T}_i\subset \{0,\ldots,T\}$;  that is, only $\{x_{i,t}\}_{t\in \mathcal{T}_i}$ is available.
In terms of the abstract dynamic norm (\ref{dynmetric3}), we handle this by leaving out those terms in the sum over $t$ in equation (\ref{dynmetric3}) that correspond to times at which data is unavailable.
Thus, the treatment of missing data we propose is not specific to fuzzy clustering.
In the fuzzy clustering framework, this corresponds to excluding those time instants $t\in\mathcal{T}_i^c$ for which $x_{i,t}$ is unavailable from both the centre update and membership likelihood update rules.
Thus, only data that is available at a particular time instant $t$ is used to calculate cluster centre coordinates at that time $t$.

To do this efficiently, we consider the known portion of trajectory $i$, namely $\{x_{i,t}\}_{t\in \mathcal{T}_i}$, as a point in the lower-dimensional space $\mathbb{R}^{d|\mathcal{T}_i|}$ for the purposes of computing Euclidean distances in the clustering algorithm.
This projection to a lower-dimensional space is easily incorporated into Algorithm 1.
For $i=1,\ldots,n$, we define $\pi_{i}:\mathbb{R}^{d(T+1)}\to \mathbb{R}^{d(T+1)}$ by $\pi_iX_i=\hat{X}_i=(\hat{x}_{i,0},\hat{x}_{i,1},\ldots,\hat{x}_{i,T})$, where
\begin{equation}
\label{hatx}
\hat{x}_{i,t}=\left\{
                \begin{array}{ll}
                  x_{i,t}, & \hbox{if $t\in \mathcal{T}_i$;} \\
                  0, & \hbox{if $t\notin\mathcal{T}_i$.}
                \end{array}
              \right.
\end{equation}

To exclude unavailable observations from centre updates, for each time instant $t=0,1,\ldots,T$, we define $\mathcal{I}_t=\{i:t\in\mathcal{T}_i\}\subset\{1,\ldots,n\}$, namely the indices of all trajectories with observations available at time $t$.

\vspace{.2cm}
\noindent\textbf{Algorithm 2:  Clustering with missing data}
\begin{enumerate}
\item Initialize membership values $u_{k,i}$.
\item Calculate centres:
\begin{equation}
\label{cupdate}
c_{k,t}=\frac{\sum_{i\in\mathcal{I}_t} u_{k,i}^m {x}_{i,t}}{\sum_{i\in\mathcal{I}_t} u_{k,i}^m},
\end{equation}
$k=1,\ldots,K$, $t=1,\ldots,T$.
Note that we take a convex combination over only those observations available at time $t$.
\item Update membership values:
\begin{equation}
\label{uupdate}
u_{k,i}=\frac{1/\|\pi_i{X}_i-\pi_iC_k\|^{2/(m-1)}}{\sum_{j=1}^K \left(1/\|\pi_i{X}_i-\pi_iC_j\|^{2/(m-1)}\right)},
\end{equation}
$k=1,\ldots,K$, $i=1,\ldots,n$.
Note that when computing Euclidean distances, we project onto only those temporal copies of phase space in which trajectory data for $X_i$ is available.
\item Evaluate the objective
 \begin{equation}
\label{objmissing}
\sum_{k=1}^K\sum_{i=1}^n u_{k,i}^m\|\pi_iX_i-\pi_iC_k\|^2.
\end{equation}
If the improvement in the objective is below a threshold, go to step 5;  otherwise go to step 2.
\item Output cluster centres $C_k\in \mathbb{R}^{dT}, k=1,\ldots,K$ and membership likelihoods $u_{k,i}\in [0,1], k=1,\ldots,K$, $i=1,\ldots,n$.
\end{enumerate}
Algorithm 2 is also frame-independent;  the proof is identical to the proof of frame-independence of Algorithm 1.

We remark that Algorithm 2 will have a preference for clusters that each contain a similar total amount of data;  for example, one cluster comprising 20 trajectories of length ten and another comprising 40 trajectories of length five both contain the same amount of data.
In some problems, one may wish Algorithm 2 to have a preference for clusters with similar numbers of trajectories, irrespective of the amount of available data in each trajectory.
To achieve this, one can replace $u_{k,i}^m$ with $u_{k,i}^m/|\mathcal{T}_i|$ in (\ref{cupdate}) in Step 2 and (\ref{objmissing}) in Step 4.
The reasoning behind this replacement is that with the factor $1/|\mathcal{T}_i|$, (\ref{objmissing}) computes the \emph{average} weighted squared distances from centres (per trajectory), whereas without this factor, the \emph{total} squared distances from centres is computed.
With this altered objective function, one constructs the correspondingly altered update rules (\ref{cupdate})--(\ref{uupdate}) as outlined below Algorithm 1.
We tested Algorithm 2 with and without this factor in the examples in Section \ref{sec:examples} and found little difference;  we report the results without this factor.

\section{What can go wrong?}\label{sec:gowrong}

Before we begin to outline some guidelines to avoid potential pitfalls in sections \ref{sec:falsepos}--\ref{sec:otherinacc}, we introduce a quantity that (along with the likelihoods $u_{k,i}$) can be useful for assessing confidence in the clustering reported by Algorithms 1 or 2.

\subsection{Entropy and classification uncertainty}

Each trajectory $x_{i,t}, t\in \mathcal{T}_i$ has relative probabilities $u_{k,i}\in [0,1]$, of belonging to cluster $C_k$, $k=1,\ldots,K$, respectively.
We can now define an overall measure of certainty of cluster assignment of trajectory $i$  via the normalised entropy of the probability vector $[u_{1,i},\ldots,u_{K,i}]$, namely
\begin{equation}
\label{eq:entropy}
h_i:=\frac{-\sum_{k=1}^K u_{k,i}\log u_{k,i}}{\log K}.
\end{equation}
The quantity $h_i$ takes values between 0 and 1, with $h_i=0$ representing certain classification of trajectory $i$ to one of the $K$ clusters and $h_i=1$ representing complete uncertainty of classification of trajectory $i$ to one of the $K$ clusters, see also Ref.~\onlinecite{bezdekbook}.

A collection of trajectories that are retained in a compact region of phase space over the time duration should correspond to a single cluster in $\mathbb{R}^{d(T+1)}$.
Each of these trajectories should therefore have a low value of $h_i$.
A spatial plot of the field $h_i$ over the phase space $\mathbb{R}^d$ is therefore useful for identifying the strength with which trajectories belong to clusters.
Finer, cluster-by-cluster spatial information can be obtained by producing $K$ spatial plots of the likelihoods $u_{k,i}$ separately for each $k=1,\ldots,K$.

\subsection{False positives}
 \label{sec:falsepos}
 Algorithms 1 and 2 will always produce centres and clusters, even if the system under consideration has no features that could be considered to be coherent.
Thus, there is the possibility of Algorithms 1 and 2 reporting false positives.
There are some easy ways to inspect the reported clusters and check for false positives.
If the phase space is in one, two, or three dimensions, then one can visually inspect the clusters at each time instant to check if the clusters do indeed mostly remain in separate compact regions.
This can be done by plotting $u_{k,i}$ against $x_{i,t}$ for $k=1,\ldots,K$ and $t=0,1,\ldots,T$ (using e.g.\ the {\tt scatter} command in MATLAB) to check the certainty of classification for individual clusters.
If the phase space is not low-dimensional, one can plot $u_{k,i}$ against $i$ (or $h_i$ against $i$) and inspect how many trajectories have high confidence of classification.
A low classification confidence is indicative of the cluster not corresponding to a coherent set.

\subsection{Choice of trajectory output times and choice of $m$}

Clustering with respect to the Euclidean metric becomes less meaningful in high dimensions, with the distribution of interpoint distances becoming increasingly tight.
This can be partly mitigated by using an $\ell_p$ norm rather than the Euclidean $\ell_2$ norm, but we have found the following rules of thumb very helpful, and have achieved good results with the standard Euclidean norm.

Firstly, one should choose the time between $x_{i,t}$ and $x_{i,t+1}$ to represent some nontrivial dynamics.
If the increment $t\to t+1$ is too short, the dynamics is close to the identity transformation, and one adds $d$ dimensions to the clustering problem (making it more difficult) for no information gain.
On the other hand, the increment from $t\to t+1$ should not be so long that the underlying dynamics appears random over one time step;  a group of nearby points at time $t$ should remain in a ``connected'' region at time $t+1$, even though this region may be stretched and folded.
Secondly, the total time duration $T$ should not be so long that the entire phase space is thoroughly mixed;  for such $T$ there is no chance of finding coherent sets.
Once the step $t\to t+1$ and total duration $T$ have been selected as above, one should obtain reasonable results.
Finally, to fine tune the value of $m$ to ensure robust results, we suggest the following rule.
Begin with $m=2$ and decrease $m$.
For each value of $m$, record the locations of $x_{i^*_k,0}$, the maximum likelihood trajectories at time $t=0$ (the choice of $t=0$ is arbitrary).
Find a range of $m$ for which the locations of the maximum likelihood trajectories are stable (i.e. approximately fixed).
Note that the centres $c_{k,0}$ at time $t=0$ will tend to continue to vary with $m$ so they are not good indicators of cluster stability with $m$.

\subsection{Centre collapse} If two or more of the reported cluster centres are all very close to one another in space, there are at least three possibilities.
Firstly, it could be that there are no coherent structures in the trajectory data.
Secondly, it could be that the choice of the step $t\to t+1$ and/or $T$ are unsuitable.
Thirdly, even if the choice of the step $t\to t+1$ and/or $T$ are reasonable, it could be that the value of $m$ is too high.
In our experiments we have found that the larger $d(T+1)$ is (the larger the total dimension), the smaller $m$ needs to be to avoid centre collapse.
This is not surprising because with higher dimension, the interpoint distances distribution is more tight, and a lower value of $m$ emphasises differences in distance more.
This is the reason behind our suggestion in the previous paragraph to start with $m=2$ and decrease $m$ until the maximum likelihood trajectories are stable.

\subsection{Other inaccurate results}
\label{sec:otherinacc}
For systems that do contain finite-time coherent sets, there are some points to bear in mind to increase the accuracy of the reported clusters.
If a finite-time coherent set is small relative to the domain size and few clusters are sought, because Algorithms 1 and 2 favour clusters containing approximately the same number of trajectories, the clusters may be much larger than the true coherent region.
In such a situation, an inspection of the likelihood functions may reveal the small coherent regions as ``high likelihood''.
On the other hand,  if there are few, large coherent sets, but one chooses a large value of $K$, then the coherent regions will likely be subdivided into several clusters.

These effects can be studied by varying the number of clusters $K$ (which is cheap to experiment with).
For each $K$ one can visually inspect the clustering confidence according to $u_{k,i}$ and $h_i$, as discussed in Section \ref{sec:falsepos}.
If a regime of cluster stability can be found for a number of consecutive $K$, this gives some confidence to the results.
Finally, if sufficient data is available, the results can be checked against the classical finite-time coherent set identification methods\cite{FSM10,F13,FPG14}.

\section{Numerical experiments}\label{sec:examples}

\subsection{One-dimensional examples}
\label{sec:1d}

We start with three one-dimensional maps on $S^1$, which we think of as the unit interval $[0,1]$ with the endpoints identified. Because we are on $S^1$ and not $[0,1]$, the distance computation and the center updating are modified in the obvious way.
The first example is given by
\begin{equation}\label{eq:3map}
S(x)=\left\lbrace \begin{array}{cc}
3x \Mod{\frac{1}{3}}+1/3, & x< \frac{1}{3}, \\
3x -\frac{1}{3}\Mod{\frac{1}{3}}+2/3, & \frac{1}{3} \leq x < \frac{2}{3},\\
3x -\frac{2}{3} \Mod{\frac{1}{3}}, &  x\geq \frac{2}{3}.
\end{array}\right.
\end{equation}
The map $S$ cyclically permutes the three intervals $[0, \frac{1}{3})$, $[\frac{1}{3}, \frac{2}{3})$ and $[\frac{2}{3}, 1)$ and mixes each interval internally.
Thus the graph of $S$ features three equally sized blocks that are cyclically permuted, see Figure \ref{fig:intervalmap}.

\begin{figure}[h!]
\begin{tabular}{cc}
\includegraphics[width=0.49\columnwidth, clip=true]{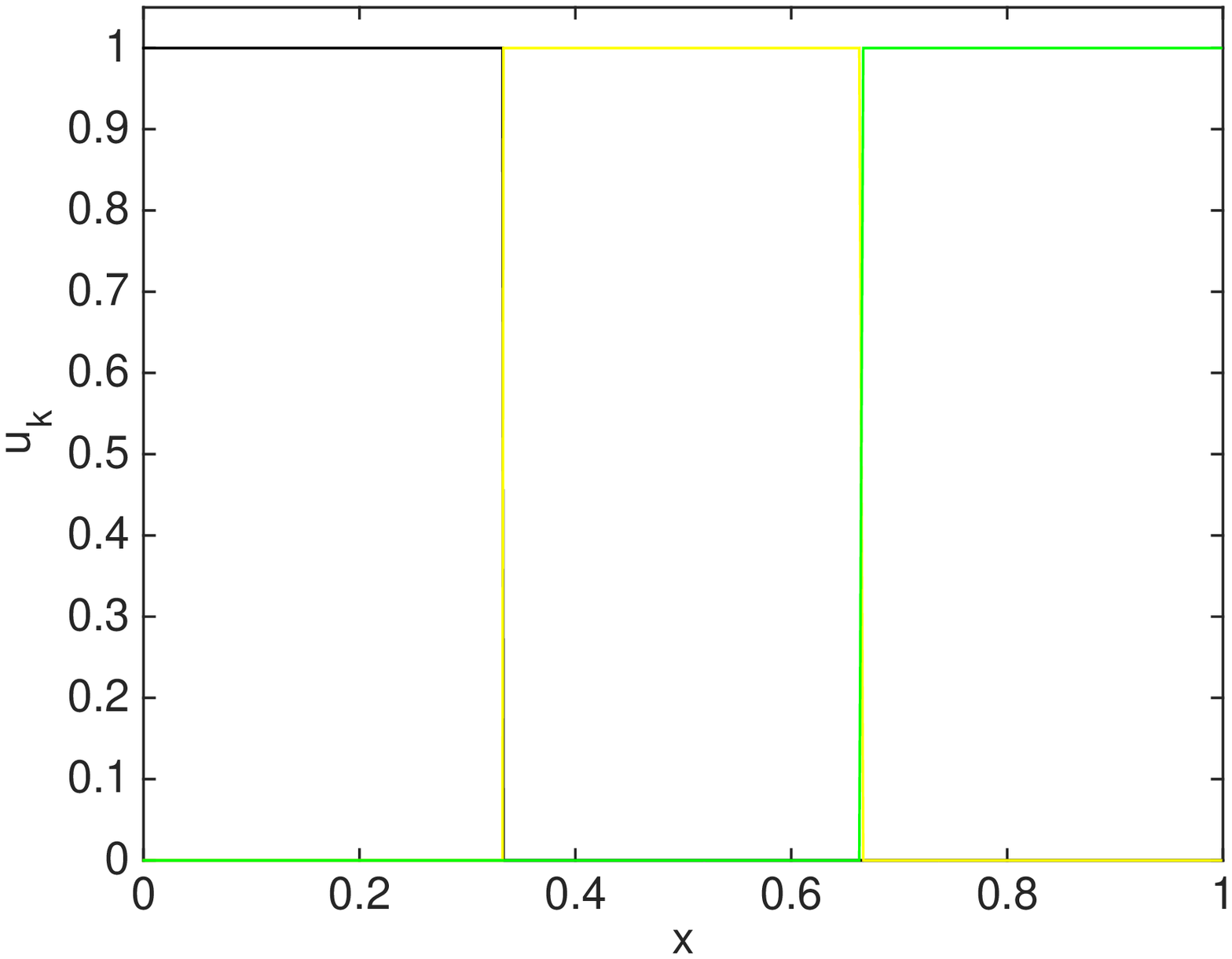} &
\includegraphics[width=0.49\columnwidth, clip=true]{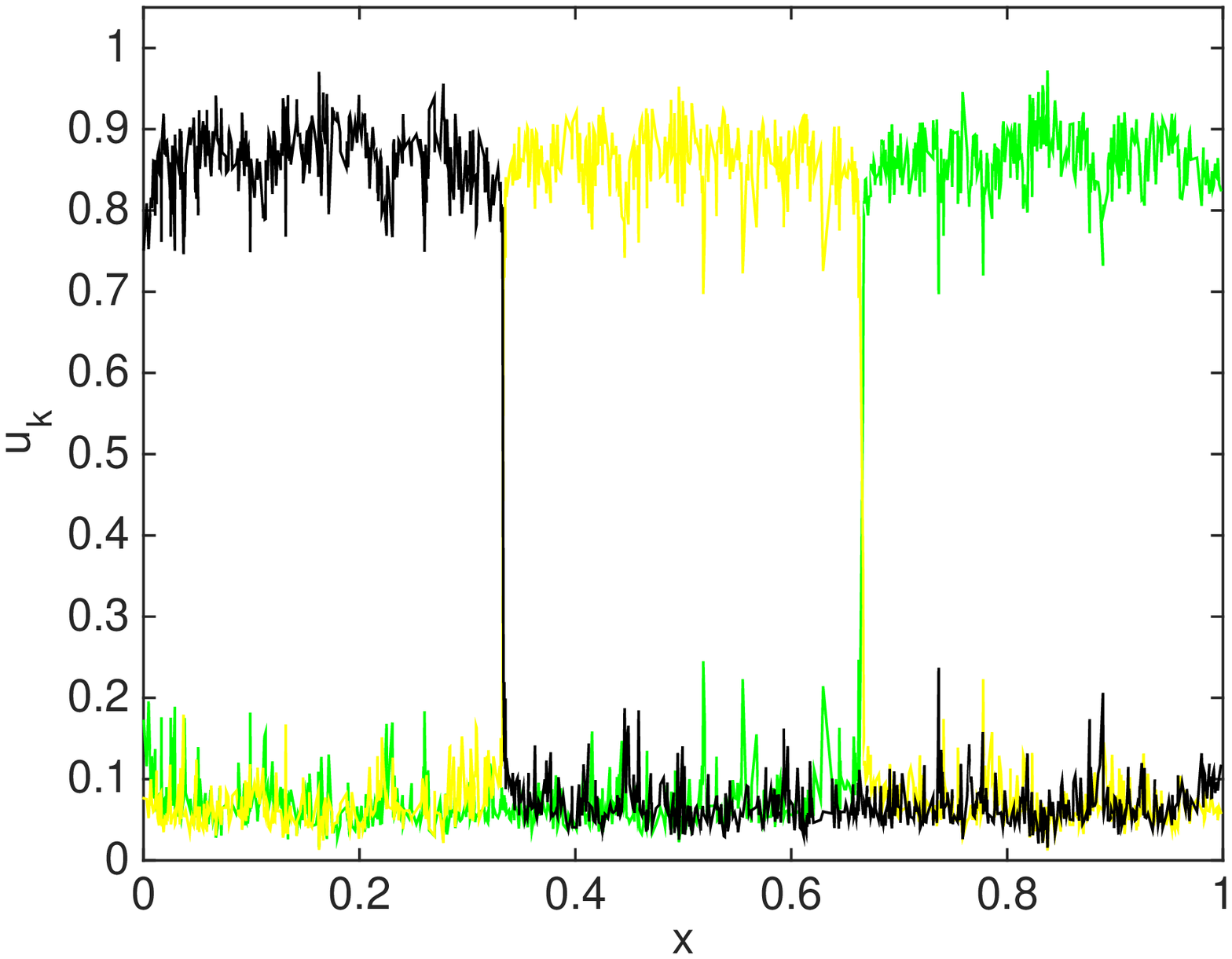}\\
 {\scriptsize (a)} & {\scriptsize (b)}\\
 \end{tabular}
\caption{Membership functions for the three clusters of the one-dimensional map (\ref{eq:3map}) plotted against the initial conditions of the trajectories under consideration.\\ (a) $m=1.1$, (b) $m=2$.}\label{fig:1dclusterm}
\end{figure}

\begin{figure}[h!]
\includegraphics[width=0.99\columnwidth, clip=true]{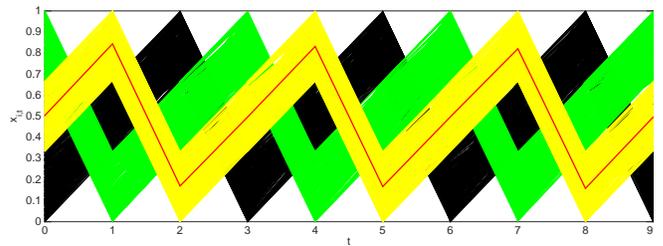}
\caption{Clustering of the 1d map (\ref{eq:3map}): time evolution of the three clusters, one of the cluster centers as computed by Algorithm 1 in red.}\label{fig:1dcluster}
\end{figure}

To test Algorithm 1 we select 1000 random initial conditions from $[0,1]$ and iterate them nine times by the mapping $S$.
We want to find clusters in 1000 data points in $\mathbb{R}^{10}$.
We choose $K=3$ and a very small fuzziness parameter of $m=1.1$.
The membership functions of the three clusters are shown in Figure \ref{fig:1dclusterm} (a).
As expected, the three coherent sets obtained are comprised of the three intervals;  the evolution of these intervals is visualized in Figure \ref{fig:1dcluster}.
The cluster centers are the centers of the intervals and the $u_{k,i}$ describe a very sharp trajectory-cluster membership.  Increasing the fuzziness to $m=2$ gives a fuzzier result, but still has clear clusters; see Figure \ref{fig:1dclusterm} (b).

\begin{figure}[h!]
\includegraphics[width=0.65\columnwidth, clip=true]{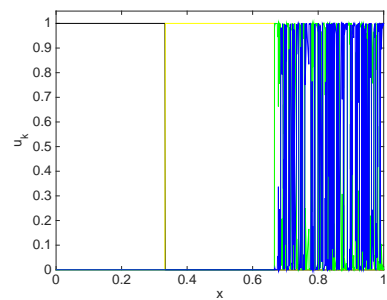}
\caption{Seeking four clusters of trajectories from (\ref{eq:3map}) results in the partition of one of the clusters in two clusters with centers almost coinciding and highly fluctuating membership functions ($m=1.1$).}\label{fig:1d4clustersentropy}
\end{figure}

If we reduce the number of desired clusters to $K=2$ the algorithm will either merge the first two or the second two clusters, depending on how the initial conditions are distributed. Trying to approximate $K=4$ coherent sets, one of the three clusters is divided into two clusters. Their centers are almost coinciding, an indication of false positives, and the membership functions on this interval are very much fluctuating, see Figure \ref{fig:1d4clustersentropy}.

\begin{figure}[h!]
\includegraphics[width=0.99\columnwidth, clip=true]{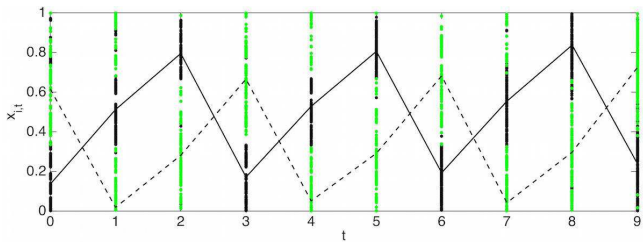}
\caption{Clustering of 1000 trajectories of length 10 of the system (FLQ10) in space-time ($m=1.5$). Plotted are only those points with $u_{k,i}>0.95$. The solid and dashed lines indicate the centers of the two clusters as computed by Algorithm 1.}\label{fig:1dclusterflq10}
\end{figure}

The map $S$ in (\ref{eq:3map}) has perfectly coherent sets:  there is no transport between the three coherent sets.
We now briefly consider two further one-dimensional systems on $S^1$ to demonstrate the more common setting of leaking coherent sets.
The first system, which will be referred to as (FLQ10), is a repeated cycle of three maps $T_1, T_2, T_3$, introduced in Ref.~\onlinecite{FLQ10}.
It was shown in Ref.~\onlinecite{FLQ10} that there are two coherent sets of different sizes that are cyclically permuted. Details of the model can be found in Ref.~\onlinecite{FLQ10} (proof of Thm. 5.1 and Figure 1).
Choosing again 1000 random initial conditions from $[0,1]$ and nine iterates of the maps (three cycles of $T_3\circ T_2\circ T_1$), we seek to find two clusters in the ten-dimensional data.
For this we choose a fuzziness constant of $m=1.5$. In Figure \ref{fig:1dclusterflq10} we show the two clusters in space-time, plotting only those points with a membership value of at least 95\% (according to the $u_{k,i}$) of belonging to one of the clusters.
As expected, the cluster centers approximately cycle with period 3. We note that the two clusters at $t=0$ (and thus at $t=3,6,9$) are consistent with the coherent sets obtained in Ref.~\onlinecite{FLQ10} (see in particular Figure 2 in  Ref.~\onlinecite{FLQ10}, where the supports of the positive/negative parts of the eigenvector shown there are in good agreement with the two clusters at $t=0$).

\begin{figure}[h!]
\includegraphics[width=0.99\columnwidth, clip=true]{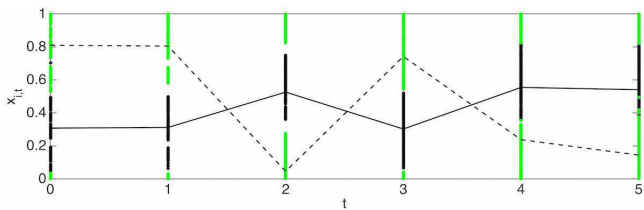}
\caption{Clustering of 1000 trajectories of the system (FLS10) in extended space ($m=1.5$).  Plotted are only those points with $u_{k,i}>0.95$. The solid and dashed lines indicate the probabilistic centers of the two clusters.}\label{fig:1dclusterfls10}
\end{figure}

\begin{figure}[h!]
\begin{tabular}{cc}
\includegraphics[width=0.49\columnwidth, clip=true]{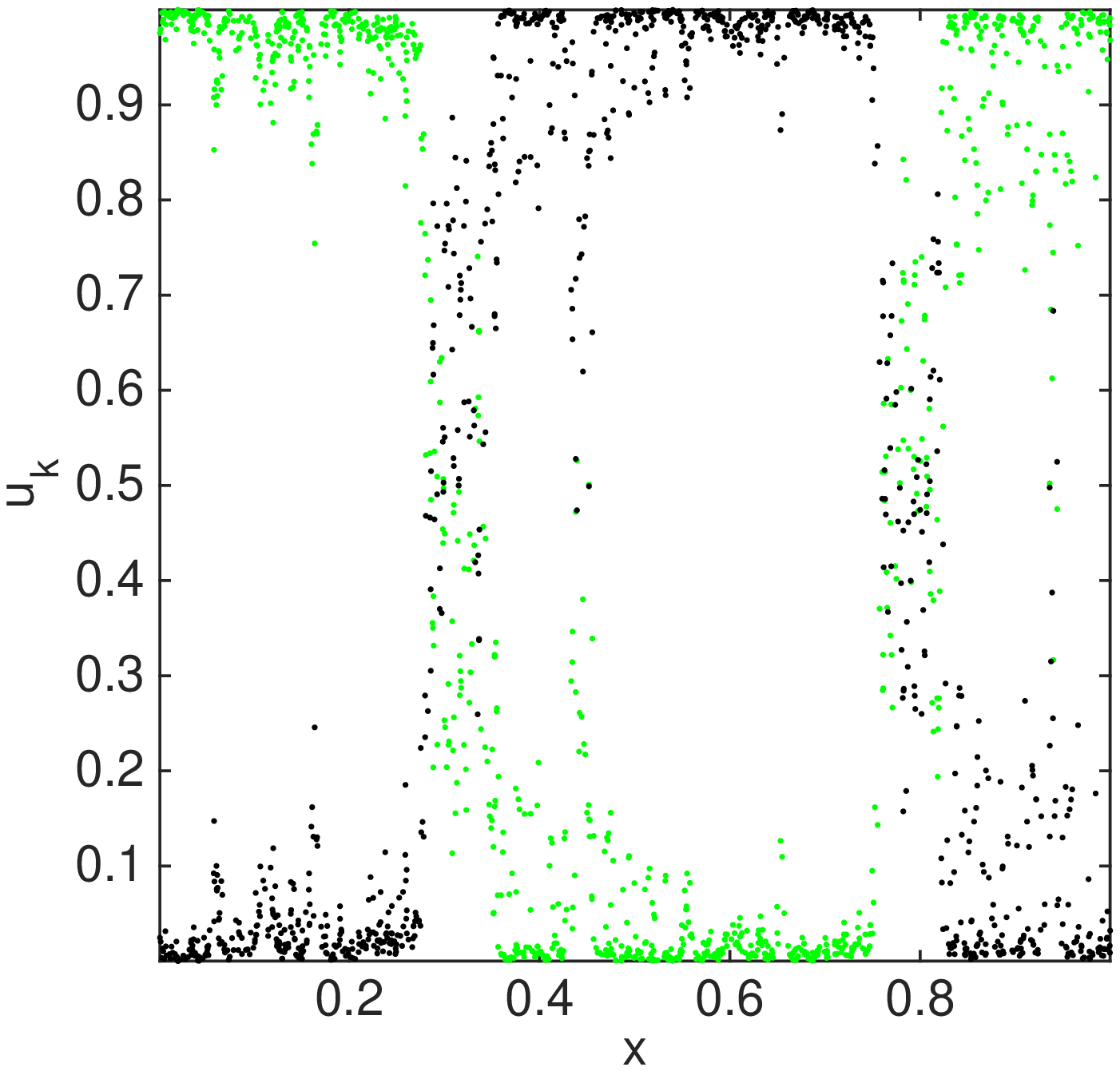}&
\includegraphics[width=0.49\columnwidth, clip=true]{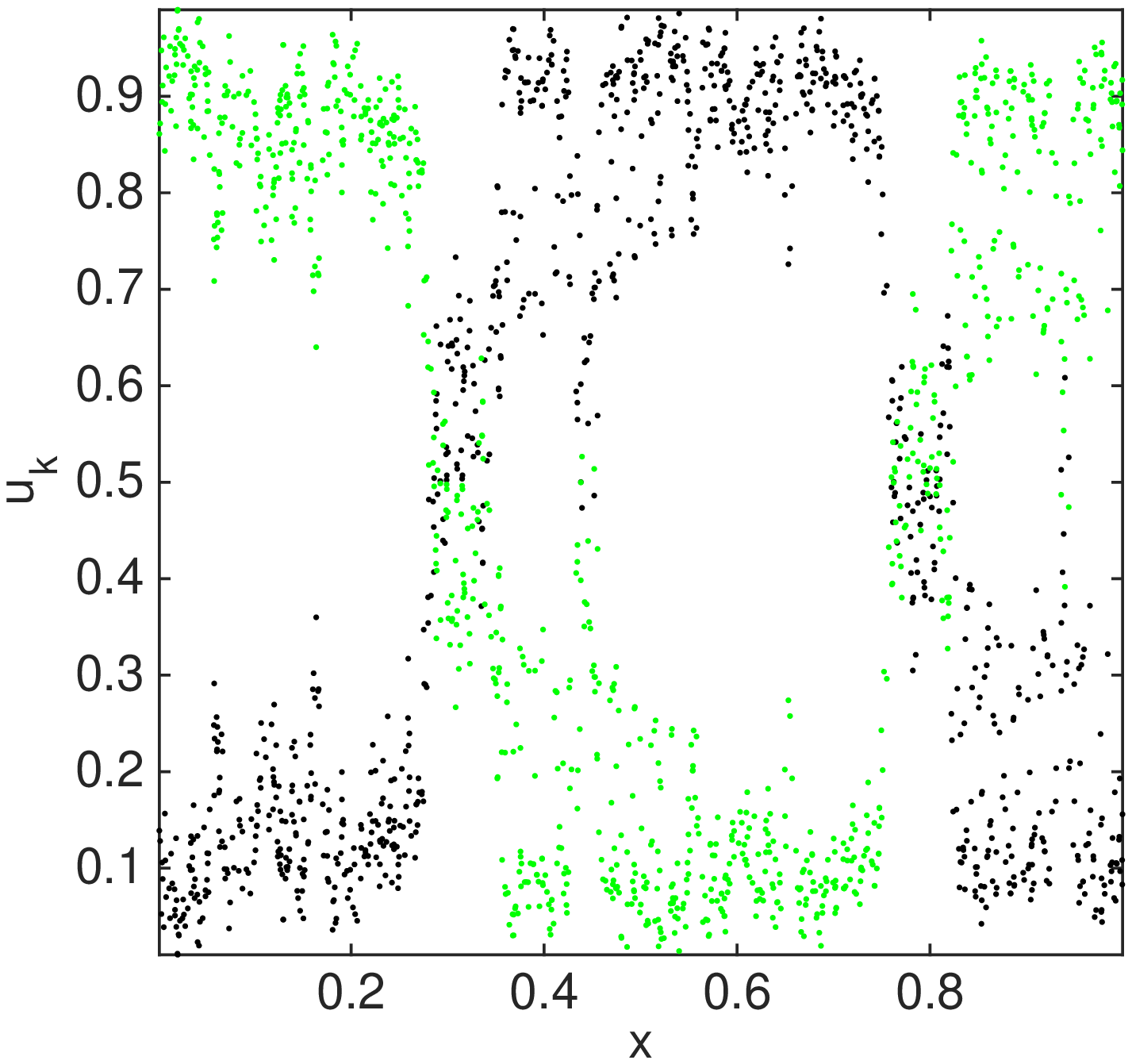}\\
 {\scriptsize (a)} & {\scriptsize (b)}\\
 \end{tabular}
\caption{Membership functions for the two clusters of  the 1d map (FLS10) at time $t=2$. (a) $m=1.5$, (b) $m=2$.}\label{fig:1dclusterfls10m}
\end{figure}

A more general situation has been discussed in Ref.~\onlinecite{FLS10}.
Two coherent sets were extracted that move in an aperiodic manner;  see Example 1 in Ref.~\onlinecite{FLS10} for more details of the underlying model, which we will refer to as (FLS10) in the following. We can reproduce the coherent behaviour of (FLS10) using the same setting as described above. Figure \ref{fig:1dclusterfls10} shows the two clusters in space-time, again plotting only those points with a membership value of at least 95\% for one of the two clusters. 
The results are consistent to those in Ref.~\onlinecite{FLS10} (see in particular Figure 8 in Ref.~\onlinecite{FLS10}, where the supports of the positive/negative parts of the Oseledets functions shown for iterates $k=0, \ldots, 5$ are in good agreement with the two clusters at times $t=0,\ldots,5$). 
The membership functions of the two clusters (plotted for time $t=2$) for the choice $m=1.5$ and $m=2$ are shown in Figure \ref{fig:1dclusterfls10m} (one can also compare the form of the black membership function with the Oseledets function in Figure 8, Ref.~\onlinecite{FLS10} for k=2).
As anticipated, the clusters are not as clear-cut as in Figure \ref{fig:1dclusterm}.
Eventually all trajectories will spread out over $[0,1]$, so that spherical compact structures as detected by our approach cease to exist.
We remark that in each of the one-dimensional examples, the maps have a uniform slope of 3, so that after the ninth iterate, nearby initial points have been separated by a factor of $3^9=19683$.

\subsection{Double gyre flow}
\label{sec:dg}

We consider the time-dependent system of differential equations \cite{shadden_lekien_marsden_05}
\begin{eqnarray}
\dot{x} &=& -\pi A \sin(\pi f(x, t))\cos(\pi y) \label{eq:dgyre}\\
\dot{y} &=& \pi A \cos(\pi f(x, t))\sin(\pi y)\frac{df}{dx}(x, t),\notag
\end{eqnarray}
where
$f(x,t)=\delta \sin(\omega t)x^2+(1-2\delta\sin(\omega t))x$.

For detailed discussions of the system we refer to Refs.~\onlinecite{shadden_lekien_marsden_05,FP09,FPG14}.
As in Refs.~\onlinecite{FP09,FPG14} we fix parameter values $A=0.25$,
$\delta=0.25$ and $\omega=2\pi$ and obtain a $1$-periodic flow. In order to be able to compare our results with those in Ref.~\onlinecite{FPG14}, where we have extracted two optimally coherent sets via transfer operator-based methods, we choose $2^{15}$ initial points on a uniform grid on the invariant set $[0, 2]\times[0,1]$. For each of these initial conditions we compute a trajectory on $[0,\tau]$, where $\tau=1,5,10$. We output the trajectory data in increments of $0.1$ time steps. Thus, for $\tau=1$ each trajectory is represented by a $22$-dimensional vector ($=(10+1) \times 2$), for $\tau=5$ and $\tau=10$ the corresponding vectors have length $102$ and $202$, respectively.

\begin{figure}[h!]
\begin{tabular}{cc}
\includegraphics[width=0.49\columnwidth, clip=true]{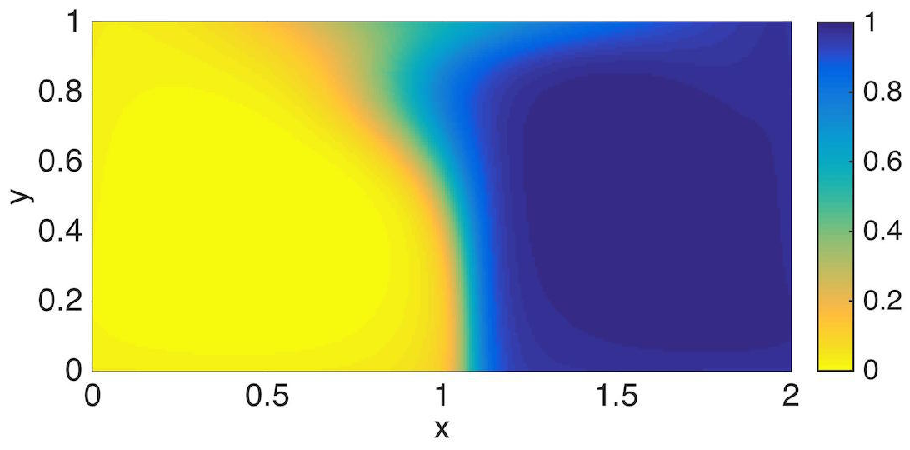} &
\includegraphics[width=0.49\columnwidth, clip=true]{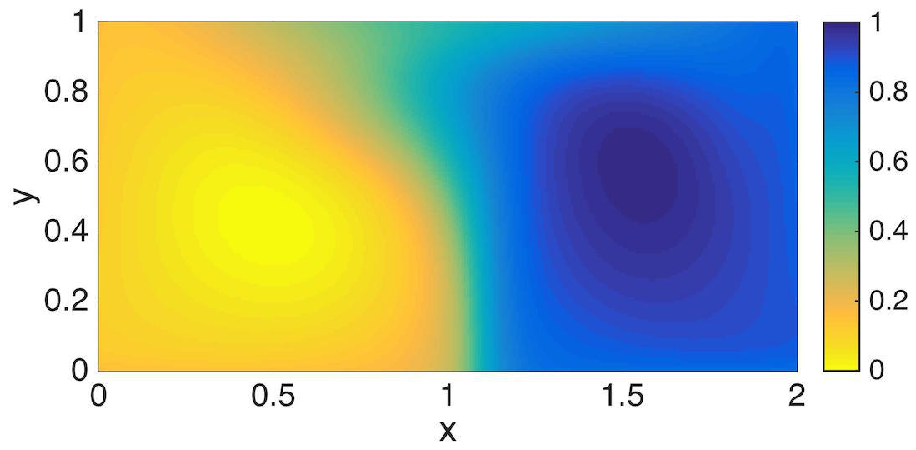}\\
 {\scriptsize (a)} & {\scriptsize (b)}\\[2mm]
\includegraphics[width=0.49\columnwidth, clip=true]{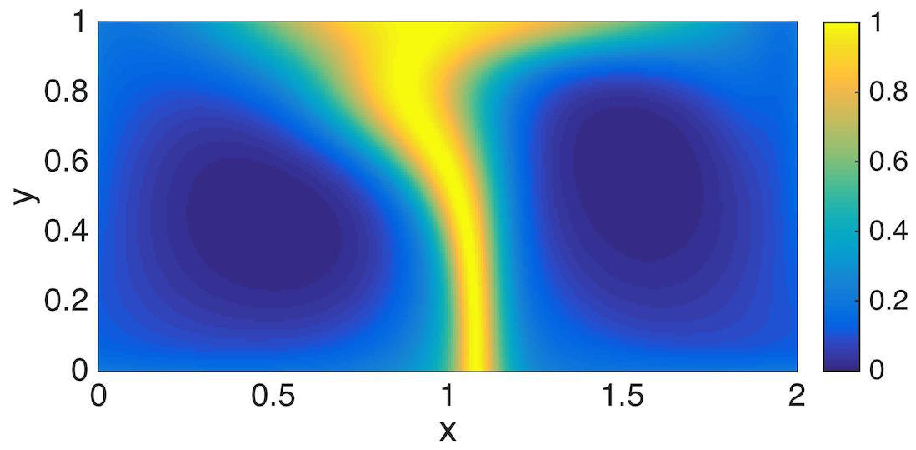} &
\includegraphics[width=0.49\columnwidth, clip=true]{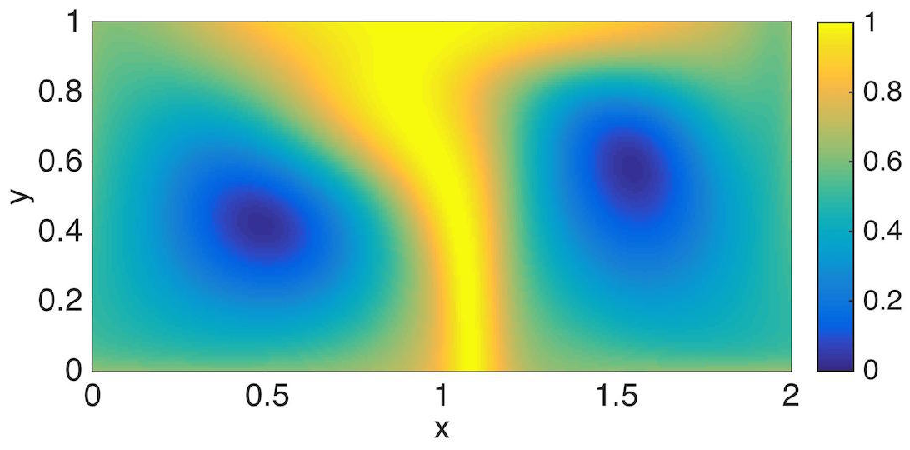}\\
 {\scriptsize (c)} & {\scriptsize (d)}\\
\end{tabular}
\caption{2-clustering of $2^{15}$ trajectories in the double gyre flow  (\ref{eq:dgyre}) with flow time $\tau=1$ for $m=1.5$ (left) and $m=2$ (right). (a),(b) membership values $u_{1,i}$ for $m=1.5$ and $m=2$, respectively. (c),(d) corresponding entropy plots using (\ref{eq:entropy}). }\label{fig:2dgyretf1}
\end{figure}

We start by extracting two clusters from the short trajectories ($\tau =1$). The upper panel of Figure \ref{fig:2dgyretf1}  displays the membership values $u_{1,i}$ (note that $u_{2,i}=1-u_{1,i}$) with respect to the initial conditions in the two dimensional phase space. To study the influence of the fuzziness exponents on the results we choose $m=1.5$ (Figure \ref{fig:2dgyretf1}(a)) and $m=2$ (Figure \ref{fig:2dgyretf1}(b)). Both plots give a clear indication of the two coherent sets. To get a more detailed picture about the certainty of cluster membership we compute the entropy $h$ from (\ref{eq:entropy}). The respective results are shown in the lower panel of Figure \ref{fig:2dgyretf1}. For the smaller fuzziness exponent $m=1.5$ (Figure \ref{fig:2dgyretf1}(c)) there are large regions of high certainty to belong to one of the two clusters, with some high uncertainty in the vicinity of the stable manifold of the hyperbolic periodic orbit on the $x$-axis. This uncertainty region increases significantly, when $m=2$ (Figure \ref{fig:2dgyretf1}(d)) is used. Here only the two regular regions (corresponding to invariant tori in the time-$1$ flow map) are highlighted as the most certain regions.

Note that a maximum likelihood hard partition into two sets gives the same result for both $m=1.5$ and $m=2$. The result is shown in Figure \ref{fig:2dgyreml} (a), with parts of the stable manifold of the hyperbolic periodic orbit on the $x$-axis superimposed. This known dominant (infinite-time) transport barrier determines a large part of the boundary between the two extracted coherent sets. This compares very well to the observations made in Ref.~\onlinecite{FPG14} (see e.g. Figure 9.3 therein).

\begin{figure}[h!]
\begin{tabular}{cc}
\includegraphics[width=0.49\columnwidth, clip=true]{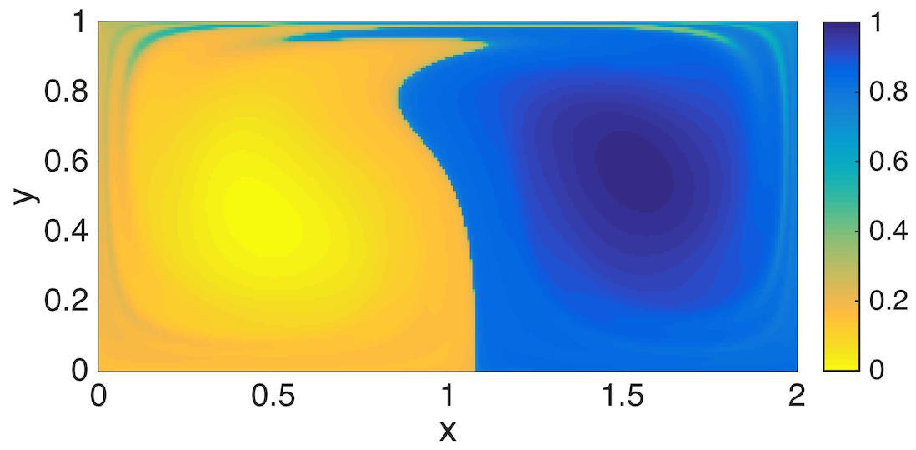} &
\includegraphics[width=0.49\columnwidth, clip=true]{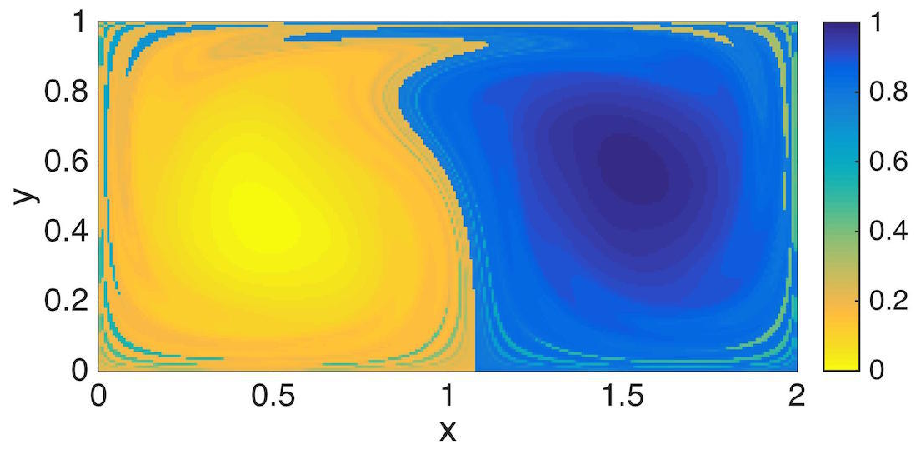}\\
 {\scriptsize (a)} & {\scriptsize (b)}\\[2mm]
\includegraphics[width=0.49\columnwidth, clip=true]{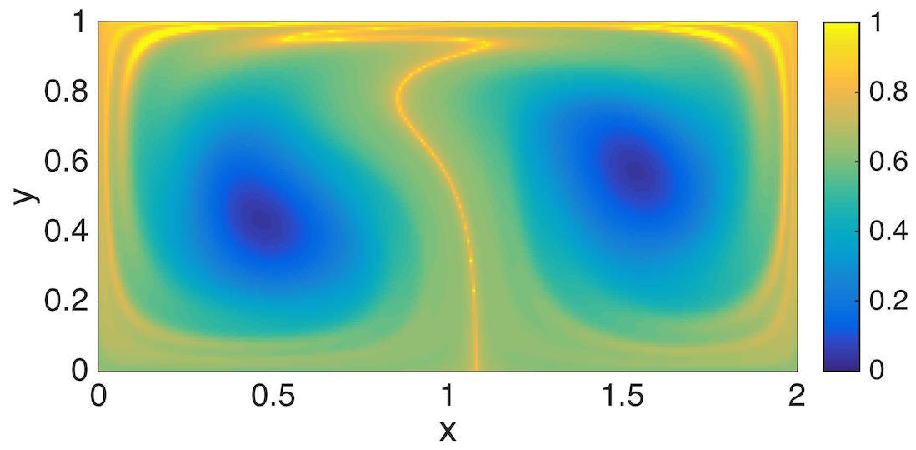} &
\includegraphics[width=0.49\columnwidth, clip=true]{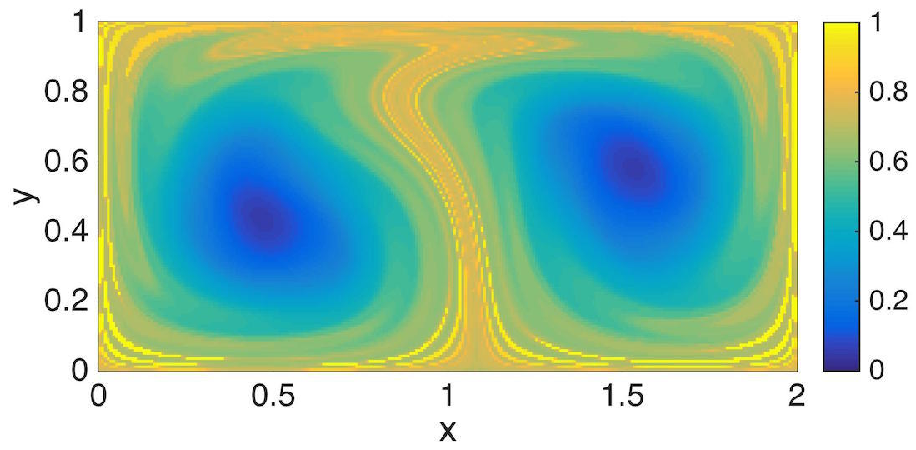}\\
 {\scriptsize (c)} & {\scriptsize (d)}\\
\end{tabular}
\caption{2-clustering of $2^{15}$ trajectories in the double gyre flow  (\ref{eq:dgyre}) for $m=2$ and flow times $\tau=5$ (left) and $\tau=10$ (right). (a),(b) membership values $u_{1,i}$ for $\tau=5$ and $\tau=10$, respectively. (c),(d) corresponding entropy plots using (\ref{eq:entropy}).}\label{fig:2dgyretf510}
\end{figure}

\begin{figure}[h!]
\begin{center}
\includegraphics[width=0.49\columnwidth, clip=true]{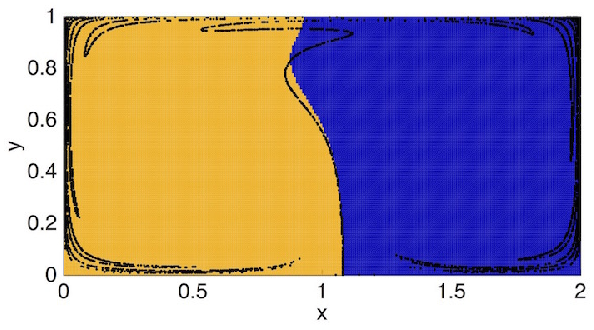}\\
{\scriptsize (a)} \\[2mm]
\begin{tabular}{cc}
\includegraphics[width=0.49\columnwidth, clip=true]{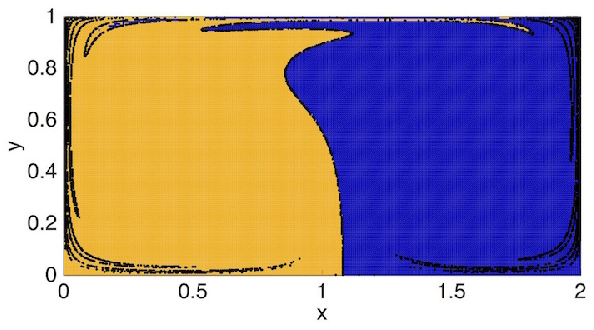} &
\includegraphics[width=0.49\columnwidth, clip=true]{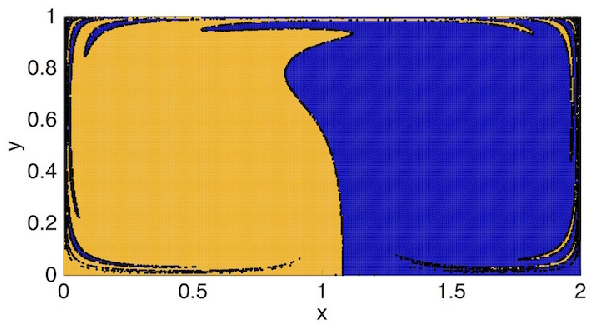}\\
{\scriptsize (b)} & {\scriptsize (c)}\\
\end{tabular}
\end{center}
\caption{Extraction of two clusters from $2^{15}$ trajectories ($m=2$) for different flow times in (\ref{eq:dgyre}) based on maximum likelihood of the membership values as in Figures \ref{fig:2dgyretf1} and \ref{fig:2dgyretf510}.  The dominant transport barrier is superimposed and bounds the two sets increasingly closely as $\tau$ increases, as demonstrated in Ref.\ \onlinecite{FPG14}. (a) flow time $\tau=1$; (b) $\tau=5$; (c) $\tau=10$.}\label{fig:2dgyreml}
\end{figure}

We now consider longer trajectories with flow times $\tau=5$ and $\tau=10$. The respective results for $m=2$ are shown in Figure \ref{fig:2dgyretf510} (a) and (b). As expected from what we have seen in Ref.~\onlinecite{FPG14} the clustering
of the initial conditions is again very much influenced by the stable manifold, see also Figure \ref{fig:2dgyreml} (b) and (c) for the respective maximum likelihood partitions into two sets.
This transport barrier also determines the regions of highest membership uncertainties, which is clearly visible in
Figure \ref{fig:2dgyretf510} (c) and (d). We note that these entropy plots have striking similarity to the finite-time entropy fields obtained directly from the transfer operator\cite{FPG12}.

A visualization of the clusters in space-time for flow time $\tau=5$ is presented in Figure \ref{fig:fc_visualize}, where from $512$ initial conditions we have plotted those trajectories for which the membership values $u_{k,i}>0.9$ ($m=2$).

\begin{figure}[h!]
\begin{tabular}{cc}
\includegraphics[width=0.49\columnwidth, clip=true]{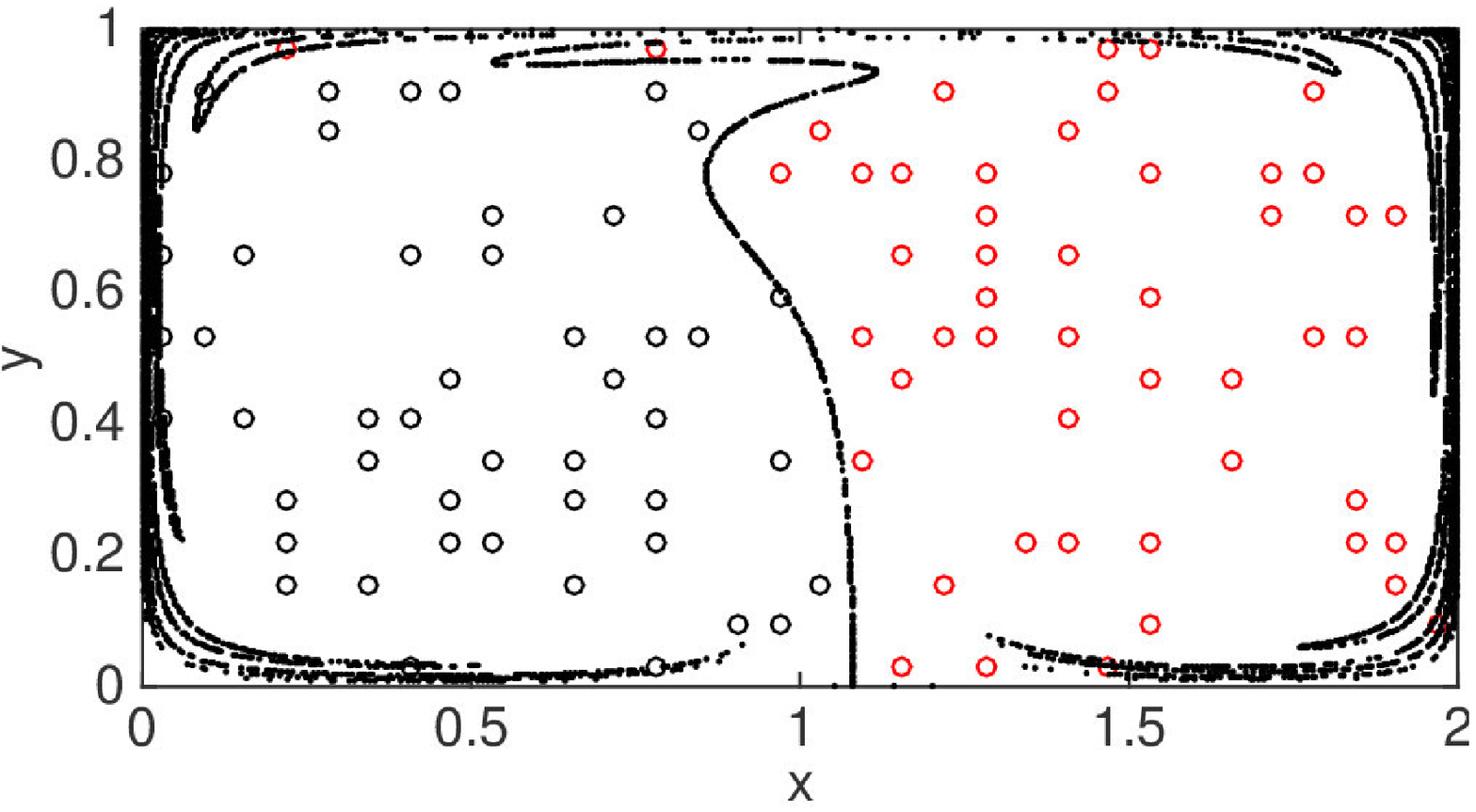} &
\includegraphics[width=0.49\columnwidth, clip=true]{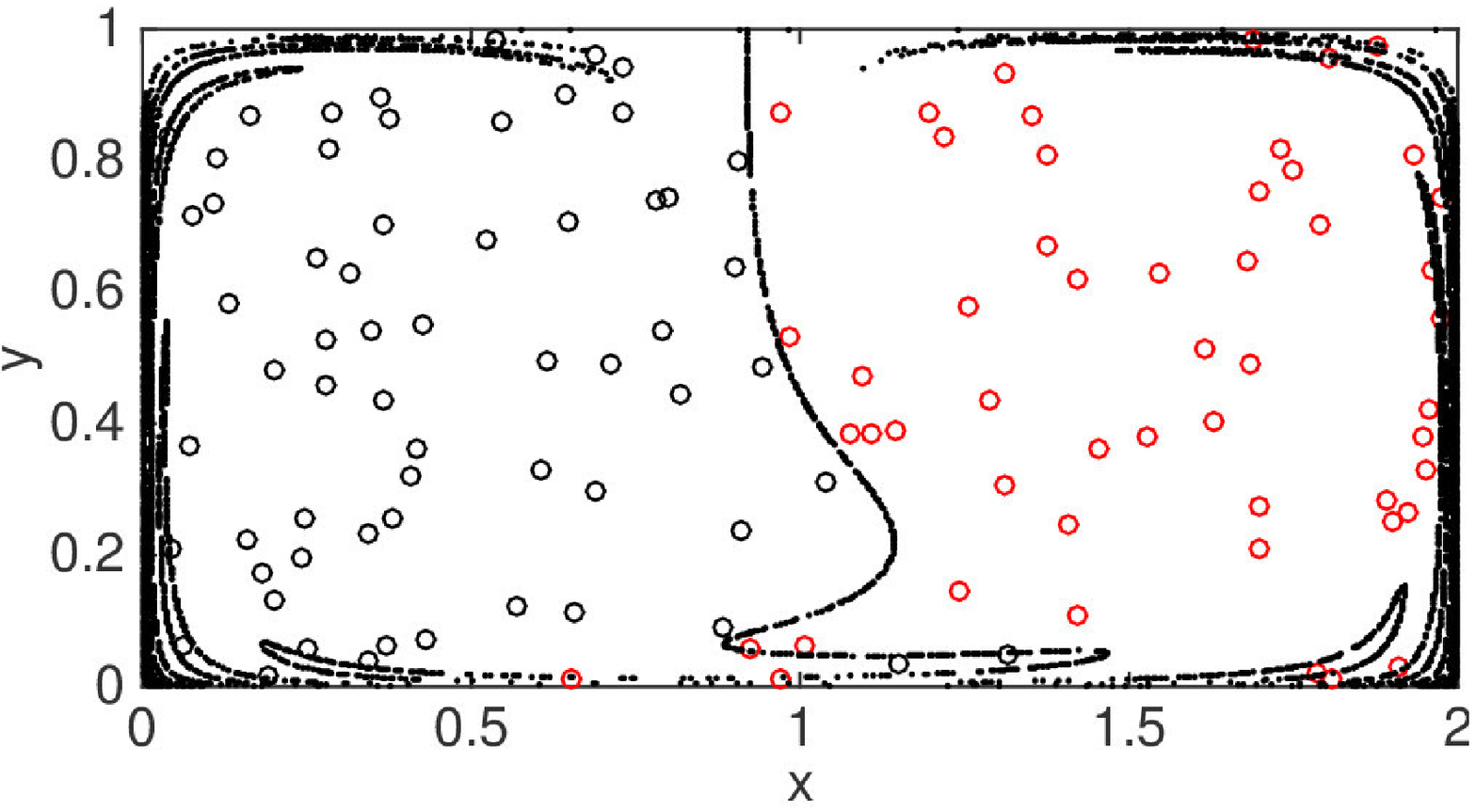}\\
{\scriptsize (a)} & {\scriptsize (b)}\\
\end{tabular}
\caption{2-clustering of double gyre flow (\ref{eq:dgyre}) with 80\% of the trajectory data missing, based on $2^9$ trajectories with flow time $\tau=5$ ($m=2$). (a) clusters at initial time, (b) clusters at final time. Corresponding transport barriers are superimposed.}\label{fig:2dgyretf1_80pc}
\end{figure}

So far we have used high-resolution and complete trajectory data. We now test our approach in the situation where the available information is poor. We use $512$ initial conditions on a regular grid on $[0, 2]\times[0,1]$ and compute trajectories for flow time $\tau=5$. We then destroy about 80\% of the trajectory information by randomly setting the particle positions to {\tt NaN}. This mimicks the situation that trajectories may not exist for the whole time span under consideration and additionally may have gaps in observation. Algorithm 2 produces two clusters from this highly incomplete trajectory data, as shown in Figure \ref{fig:2dgyretf1_80pc}.
Note that even with this severe data thinning, Algorithm 2 still classifies the remaining data points to the correct sides of the transport barriers.

\begin{figure}[h!]
\begin{tabular}{cc}
\includegraphics[width=0.49\columnwidth, clip=true]{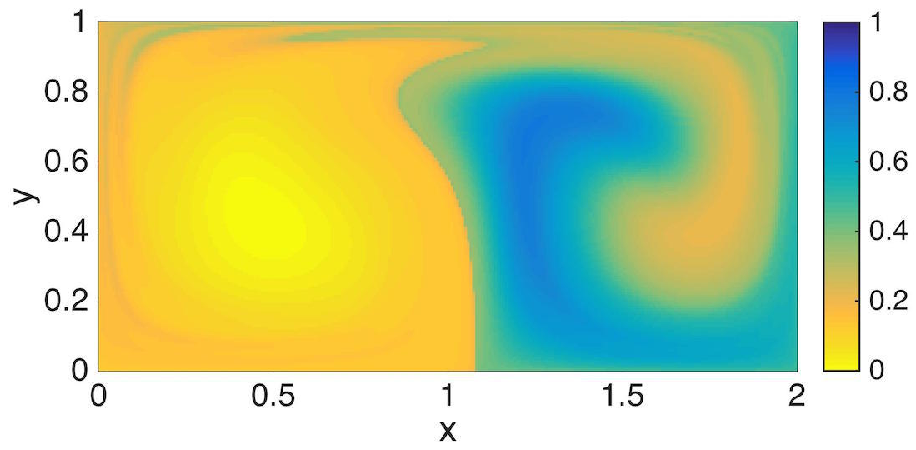} &
\includegraphics[width=0.49\columnwidth, clip=true]{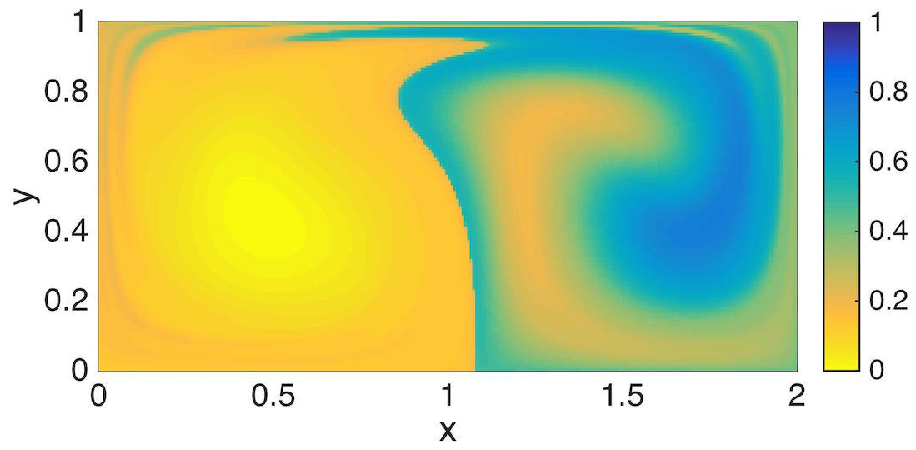}\\
{\scriptsize (a)} & {\scriptsize (b)}\\[2mm]
\end{tabular}
\includegraphics[width=0.49\columnwidth, clip=true]{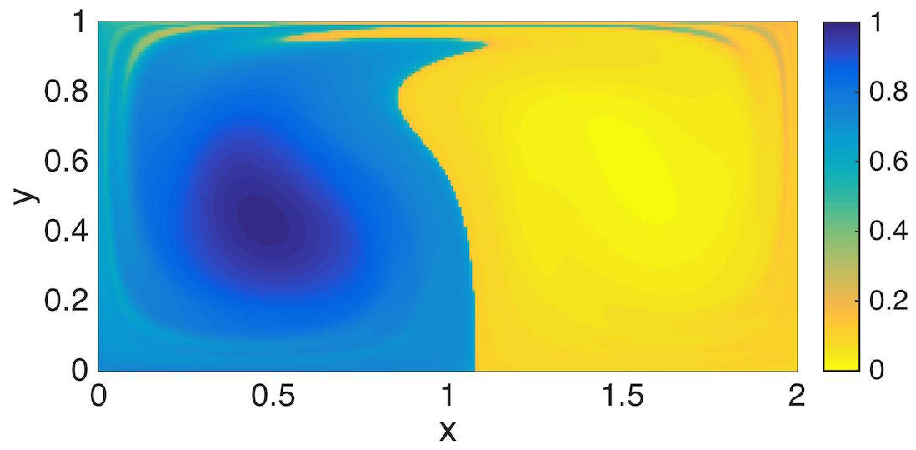} \\
{\scriptsize (c)}
\caption{3-clustering of $2^{15}$ trajectories in the double gyre flow for flow time $\tau=5$ and $m=2$. (a-c) membership functions for the three clusters.}\label{fig:2dgyreentrop5}
\end{figure}

Finally, we test what happens if we set $K>2$. We restrict again to flow time $\tau=5$ and $2^{15}$ trajectories and choose $m=2$. If $K=3$ then compared to $K=2$ either the left or right cluster is subdivided as seen in the membership values in Figure \ref{fig:2dgyreentrop5} (a-c).
The shapes of the resulting clusters in Figure \ref{fig:2dgyreentrop5} (a,b) do not have any similarity with known coherent structures for this system, but apparently the respective trajectory bundles stay coherent in our sense - with the cluster centers well separated.
However, Figure \ref{fig:2dgyreentrop5} (c) reveals that the left cluster, which is also present in the 2-clustering considered in Figure \ref{fig:2dgyretf510}, is characterized by much higher membership function values compared to the other two clusters.
We note that for $K=4$ we get a similar picture with the former two clusters both divided into two parts, and for $K=5$ one of the former two clusters is divided into two and the other into three parts.

\subsection{Transitory double gyre flow}
We consider the transitory dynamical system\cite{mosovsky_meiss_2011}
\begin{equation} \label{eqn:tdgyre}
\dot{x} = -\frac{\partial}{\partial y} \Psi,\quad \dot{y} = \frac{\partial}{\partial x} \Psi
\end{equation}
with stream function
\begin{eqnarray*}
\Psi(x,y,t)&=&(1-s(t))\Psi_P + s(t)\Psi_F\\
\Psi_P(x,y) &=& \sin(2\pi x) \sin(\pi y) \\
\Psi_F(x,y)&=& \sin(\pi x) \sin(2\pi y)
\end{eqnarray*}
and transition function
$$
s(t)=\left\lbrace \begin{array}{cc}
0, & t<0, \\
t^2(3-2t), & 0 \leq t \leq 1,\\
1, &  t>1.
\end{array}\right.
$$
The nonautonomous dynamics, which rotates a double gyre pattern counter-clockwise by 90 degrees, is restricted to the time-interval $[0,1]$. In Ref.~\onlinecite{FPG14} we have analysed this system using the transfer operator based coherent set framework. The unit square  $[0,1]^2$ is invariant under the flow and we choose $2^{14}$ initial conditions on a regular grid. We consider the flow on the transition interval $[0,1]$ and output the trajectory data in increments of $0.1$ time steps. So $T=10$ and thus we represent every trajectory as a $22$-dimensional vector.

\begin{figure}[h!]
\begin{tabular}{cc}
\includegraphics[width=0.49\columnwidth, clip=true]{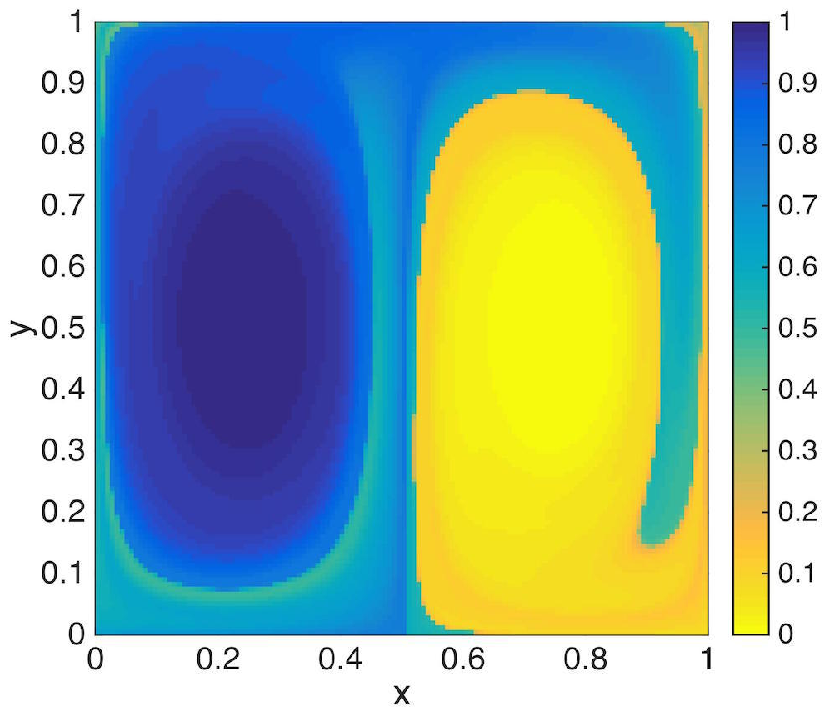} &
\includegraphics[width=0.49\columnwidth, clip=true]{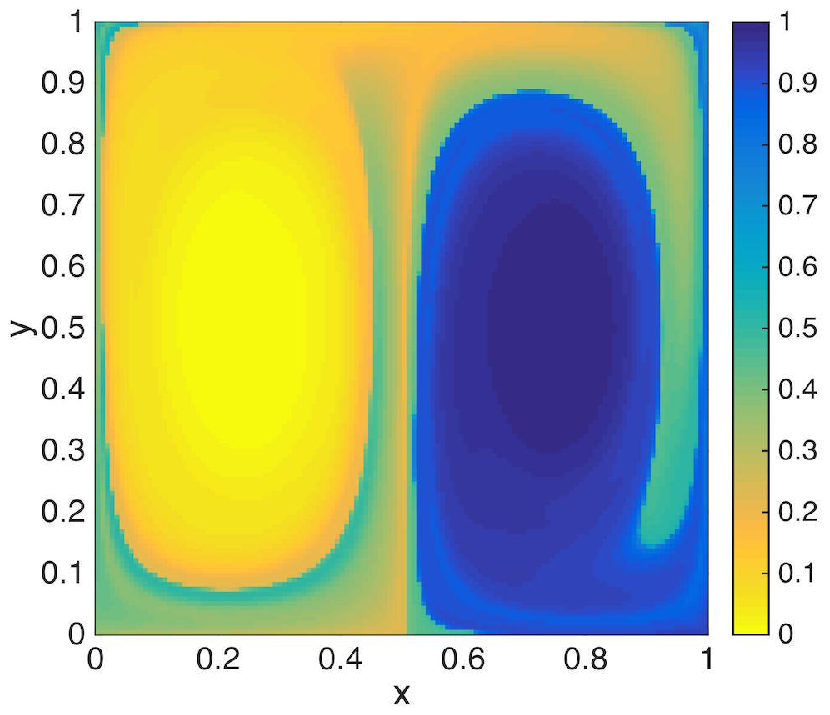}\\
{\scriptsize (a)} & {\scriptsize (b)}\\[2mm]
\end{tabular}
\caption{2-clustering of $2^{14}$ trajectories in the transitory double gyre flow (\ref{eqn:tdgyre}) for the time interval  $[0,1]$.  (a),(b) membership values $u_{k,i}$, $k=1,2$ for $m=1.5$.}\label{fig:rotdg_2clusters}
\end{figure}

Algorithm 1 with $K=2$ and $m=1.5$ returns $u_{k,i}$ that take high values on the coherent sets observed in Ref.~\onlinecite{FPG14}, as shown in Figure \ref{fig:rotdg_2clusters}.

\begin{figure}[h!]
\includegraphics[width=0.75\columnwidth, clip=true]{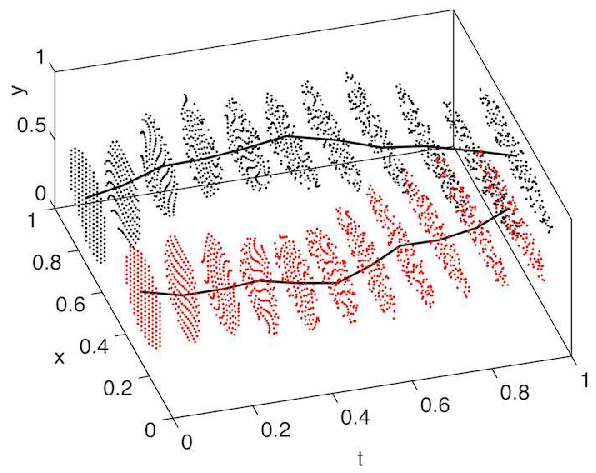}
\caption{Space-time plot of the two clusters obtained from $2^{10}$ trajectories in the transitory double gyre flow. Only trajectories with $u_{k,i}>0.95$ are shown, where $m=1.5$. The solid curves indicate the position of the probabilistic cluster centers.}\label{fig:rotdg_3d}
\end{figure}

%A partition of phase space into two clusters according to maximum likelihood is shown in Figure \ref{fig:rotdg_ml}  (a) - the results obtained for $m=1.5$ and $m=2$ are indistinguishable.
A visualization of the two clusters in space-time is presented in Figure \ref{fig:rotdg_3d}, where from $1024$ initial conditions we have plotted those trajectories for which $u_{k,i}>0.95$  ($m=1.5$) together with the probabilistic cluster centers.

\begin{figure}[h!]
\begin{tabular}{cc}
\includegraphics[width=0.49\columnwidth, clip=true]{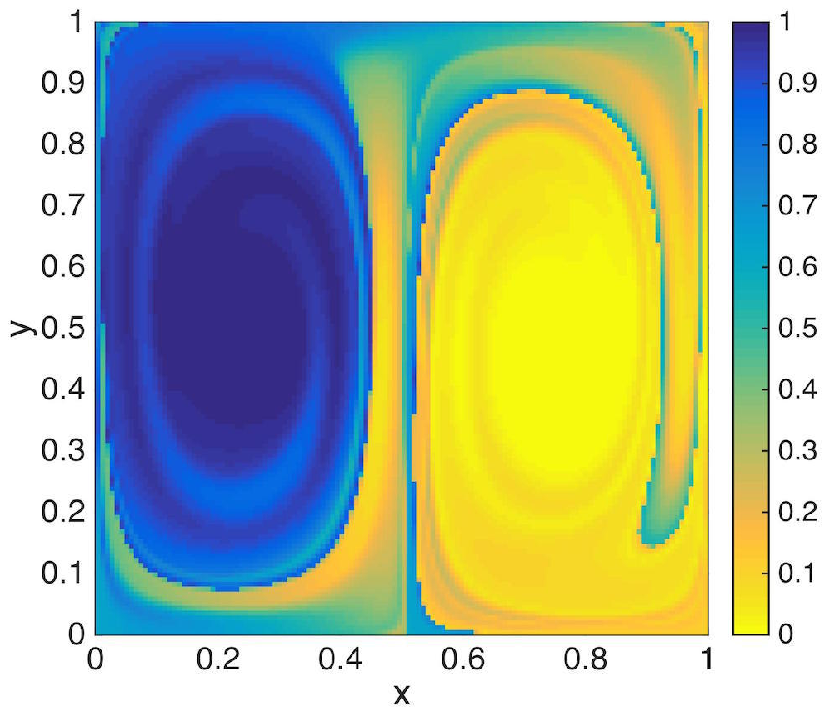} &
\includegraphics[width=0.49\columnwidth, clip=true]{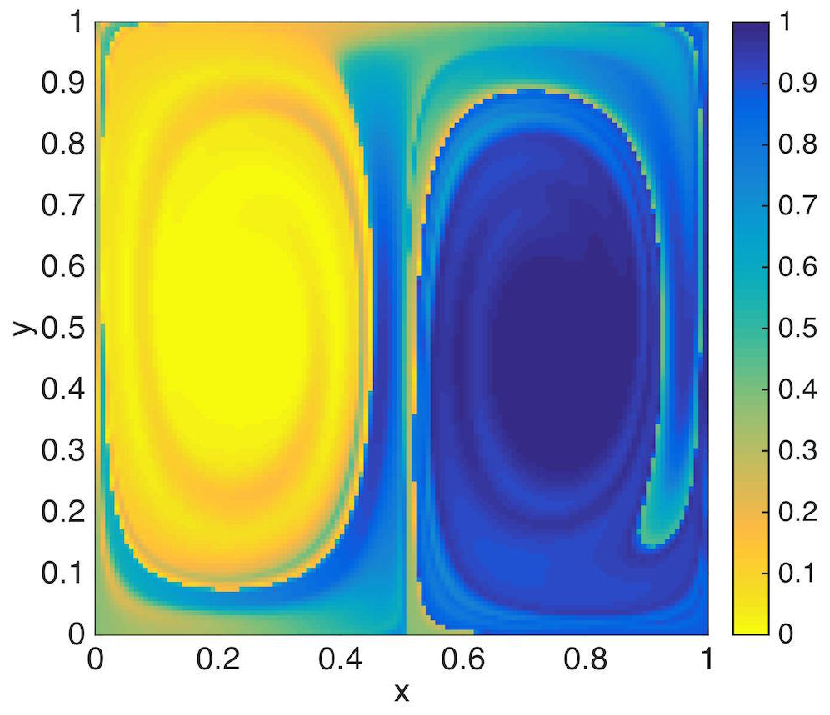}\\
{\scriptsize (a)} & {\scriptsize (b)}\\
\end{tabular}
\caption{Membership functions for the 2-clustering of transitory double gyre flow -- using only initial and end points of the trajectories ($m=1.5$).}\label{fig:rotdg_endonly}
\end{figure}

We study the influence of using information along a trajectory instead of only considering the initial and final points of a trajectory as many other identification algorithms do (i.e. taking $11$ vs $2$ time instances on $[0,1]$). The results of clustering  $2^{14}$ trajectories based only on the initial and end points of the trajectories are shown in Figure \ref{fig:rotdg_endonly}. The clusters are less smooth;  an intuitive explanation for this is that Algorithm 1 only uses point information, as opposed to probability flow information as in Refs.~\onlinecite{FSM10,F13,FPG14}.
Algorithm 1 needs to compensate for this by augmenting the point information with additional points over time.

%which is in particular highlighted in the maximum likelihood partition in Figure \ref{fig:rotdg_ml} (b). Note that like in Ref.~\onlinecite{FPG14} the left blob is picked up, whereas for $T=11$ it is the right blob (Figure \ref{fig:rotdg_ml} (a)).
%

%\begin{figure}[h!]
%\begin{tabular}{cc}
%\includegraphics[width=0.42\columnwidth]{fig17a}  &
%\includegraphics[width=0.42\columnwidth]{fig17b}\\
%{\scriptsize (a)} & {\scriptsize (b)}\\
%\end{tabular}
%\caption{Extraction of two clusters from $2^{14}$ trajectories on $[0,1]$ ($m=2$) in the transitory double gyre flow based on maximum likelihood of the membership values. (a) $T=11$ as in Figure \ref{fig:rotdg_2clusters}; (b) $T=2$ (i.e. using only initial and final points of the trajectories) as in Figure \ref{fig:rotdg_endonly}.}\label{fig:rotdg_ml}
%\end{figure}

%As each of the high certainty parts (in terms of the $u_{k,i}$) of the extracted coherent sets covers about $1/3$ of phase space we attempt to obtain a more clear-cut picture by seeking to find three clusters. This turns out to be infeasible with our approach: the cluster centers remain approximately the same as for $K=2$ but now with two of the three centers almost exactly coinciding. So in our cluster-based framework the transitory double gyre flow exhibits exactly two regions that remain compact under evolution.

\subsection{Drifter data}
We demonstrate the efficacy of our approach on real-world data, namely drifter data from the Global Ocean Drifter Program available at AOML/NOAA Drifter Data Assembly Center ({\tt http://www.aoml.noaa.gov/envids/gld/}).
The entire dataset spans the years 1979--2014, with drifter positions given every six hours. The area of observation is the global ocean (latitude $[90,-78]$ and longitude $[-180,180]$). We focus on the years 2005--2009 and restrict to those drifters that have a minimum lifetime of one year within this five-year time span.
We output the position of these 2267 trajectories (in longitude, latitude coordinates) every month, i.e. the length of our trajectories is $60$ months.

\begin{figure}[h!]
\begin{tabular}{cc}
\includegraphics[width=0.49\columnwidth, clip=true]{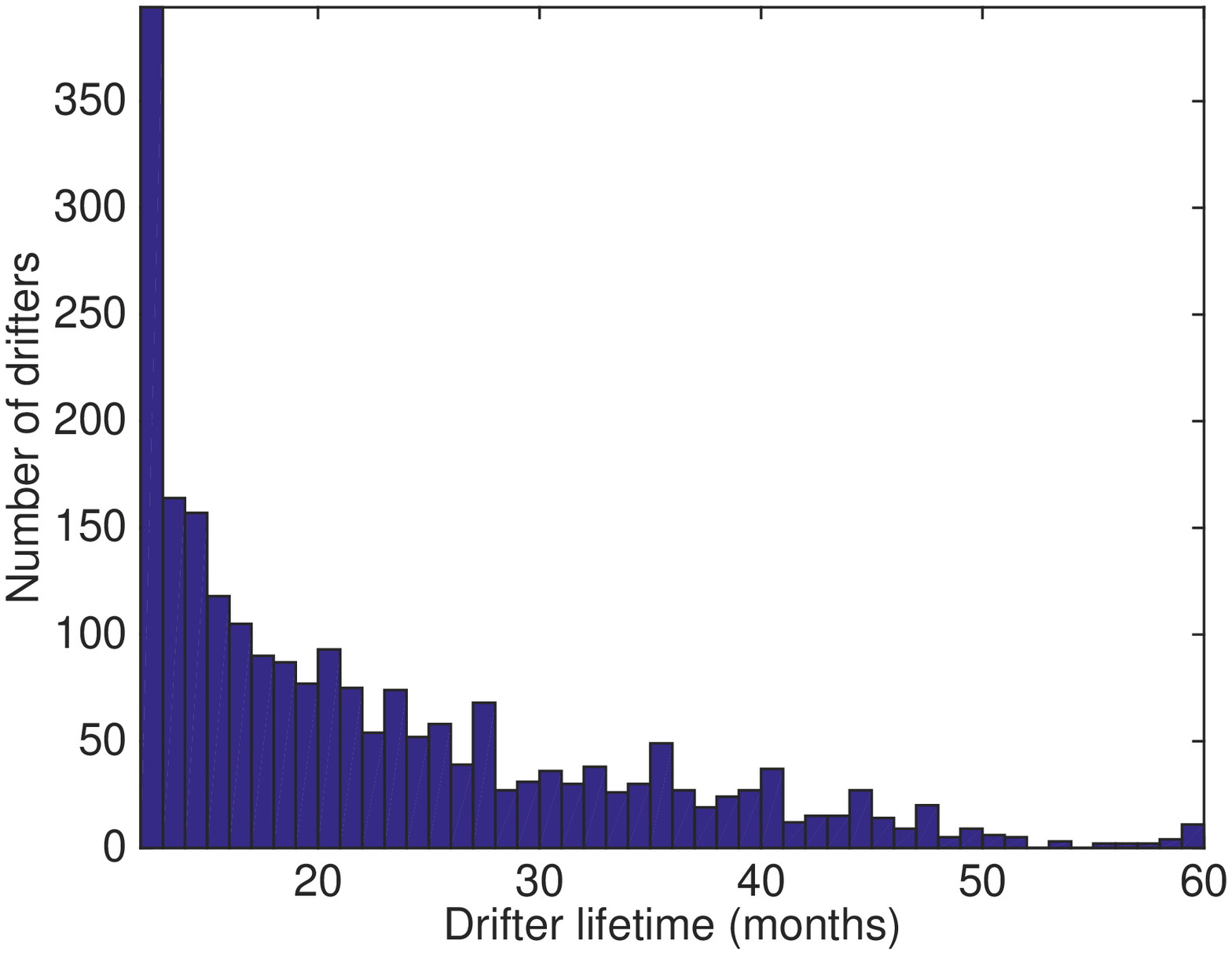} &
\includegraphics[width=0.49\columnwidth, clip=true]{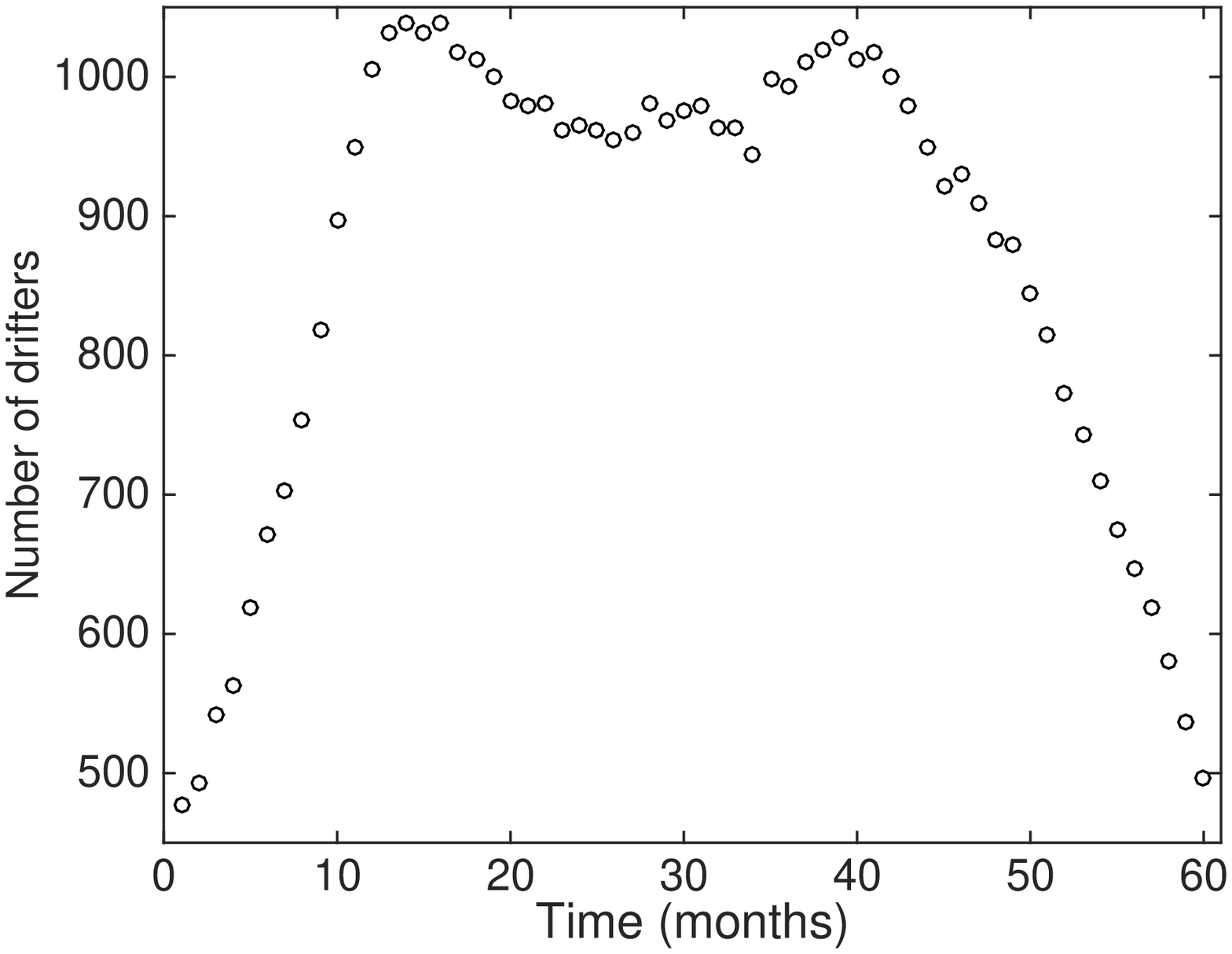}\\
{\scriptsize (a)} & {\scriptsize (b)}\\
\end{tabular}
\caption{Drifter statistics. (a) histogram of the drifter lifetimes, (b) number of drifters available at a certain time instance.}\label{fig:drifterstats}
\end{figure}

We note that a typical drifter does not operate over the whole five years;  that is, many terminate prior to December 2009 and many begin later than January 2005.
There are also gaps in observations when there is a failure in recording the drifter location, so the data is highly incomplete.
Figure \ref{fig:drifterstats} summarises two statistics:  the distribution of drifter lifetimes and the number of drifters actively recording each month.
The average lifetime of a drifter in this data set is about 23 months, with many drifters operating only for a year and only very few drifters for 4--5 years, see Figure \ref{fig:drifterstats} (a).
On average, 869 trajectories (or 38\% of all drifters in the period 2005--2009) are available at a given time instant, with less data at the beginning and the end of the considered five year time span; see Figure \ref{fig:drifterstats} (b).

As we consider the global ocean we have to respect distances on a sphere (we assume the surface of the ocean to be approximately spherical).
We also have to ensure that we restrict cluster centers to the surface of this sphere.
To achieve both of these requirements we use a cosine distance function, and update centers only on the sphere\cite{dhillon_modha_01}.
Every drifter trajectory is represented as a vector in $3 \times 60 =180$-dimensional space.
In contrast to our calculations in Algorithm 2, we simply display our results in cartesian longitude-latitude coordinates.

\begin{figure}[h!]
\begin{tabular}{c}
\includegraphics[width=0.98\columnwidth, clip=true]{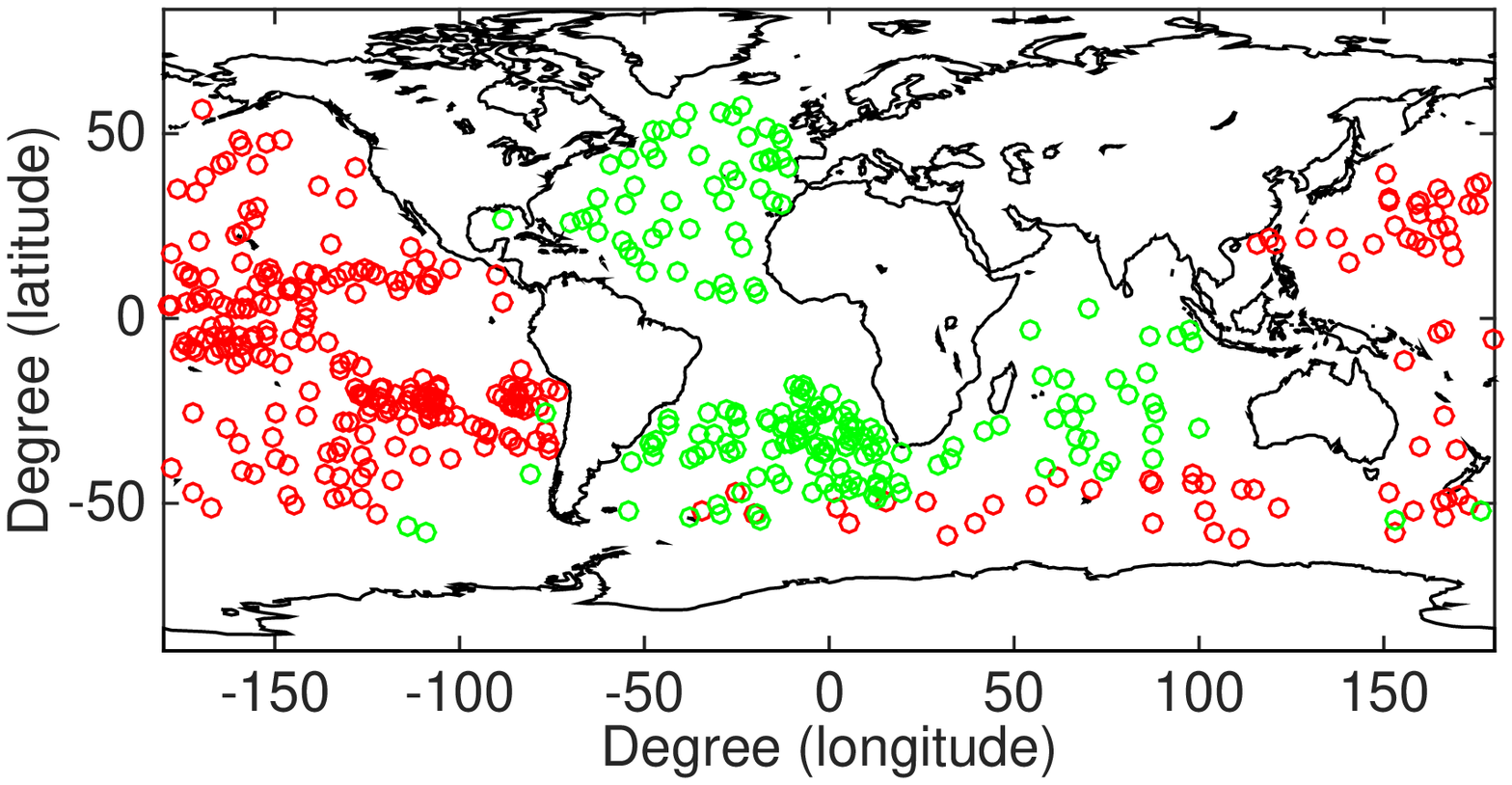} \\
{\scriptsize (a)}\\[2mm]
\includegraphics[width=0.98\columnwidth, clip=true]{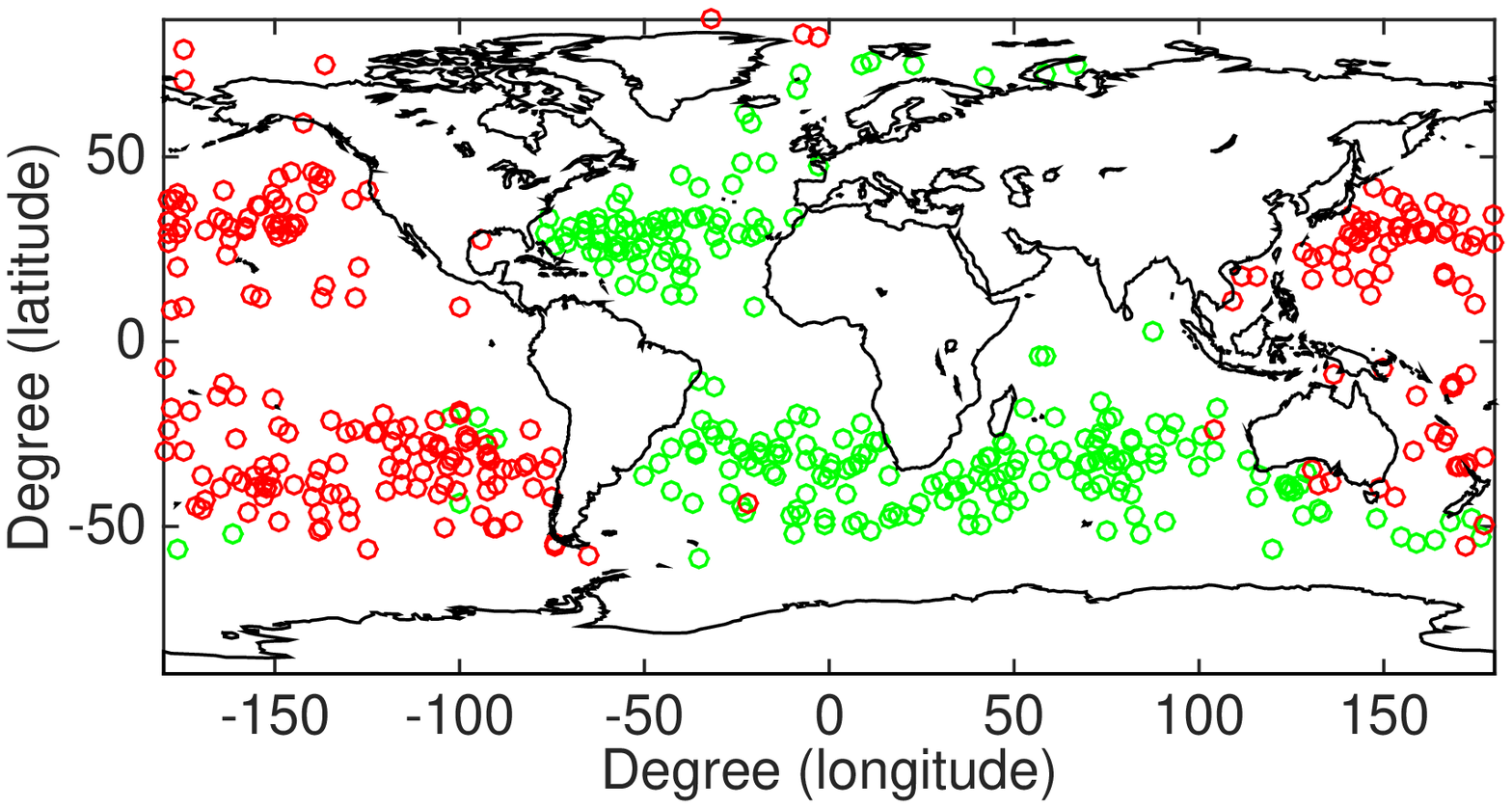}\\
{\scriptsize (b)}\\
\end{tabular}
\caption{2-clustering of drifter data ($m=1.5$).\\ (a) drifter positions January 2005, (b) drifter positions December 2009. Animation of the cluster motion 2005--2009 in online version (Multimedia view).}\label{fig:drifter2}
\end{figure}

We first look for two clusters;  Figure \ref{fig:drifter2} shows results of the clustering algorithm for $K=2$.
Figure \ref{fig:drifter2} (a) shows all drifter positions available on January 2005, coloured according to their maximum likelihood membership in one of the two clusters.
Figure \ref{fig:drifter2} (b) shows all drifter positions available on December 2009, again coloured according to their most likely cluster membership.
Thus, we expect the red (resp.\ green) cluster in Figure \ref{fig:drifter2} (a) to evolve coherently to the red cluster in Figure \ref{fig:drifter2} (b).
Of course, many of the drifters in Figure \ref{fig:drifter2} (a) do not correspond to the same physical drifter in Figure \ref{fig:drifter2} (b) because the lifetimes of many drifters are shorter than five years.
Nevertheless, as physical drifters enter and leave the dataset over the five-year duration, the drifters tagged red (resp.\ green) move as a coherent cloud.
This is illustrated in a video attached to our electronic submission, see Figure \ref{fig:drifter2} (Multimedia view).

In Figure \ref{fig:drifter2}, one sees a separation of the Pacific Ocean (red) from the Atlantic and Indian Oceans (green), which are grouped together.
Here, continental obstructions play an obvious role in the dynamical separation of the ocean surface flow.
Figure \ref{fig:drifter2} (a) ascribes the southern part of the Indian Ocean to the Pacific Ocean.
This is in line with recent research\cite{FSvS14} (see Figure 6 in Ref.~\onlinecite{FSvS14}) based on transfer operator analysis of the Ocean General Circulation Model for the Earth Simulator (OFES model)\cite{Masumoto04,Sasaki08}, and consistent with a general eastward flow of water in high southern latitudes.
One observes that the red drifters in Figure \ref{fig:drifter2} (a) have flowed eastwards to rejoin the Pacific in Figure \ref{fig:drifter2} (b).

\begin{figure}[h!]
\begin{tabular}{c}
\includegraphics[width=0.98\columnwidth, clip=true]{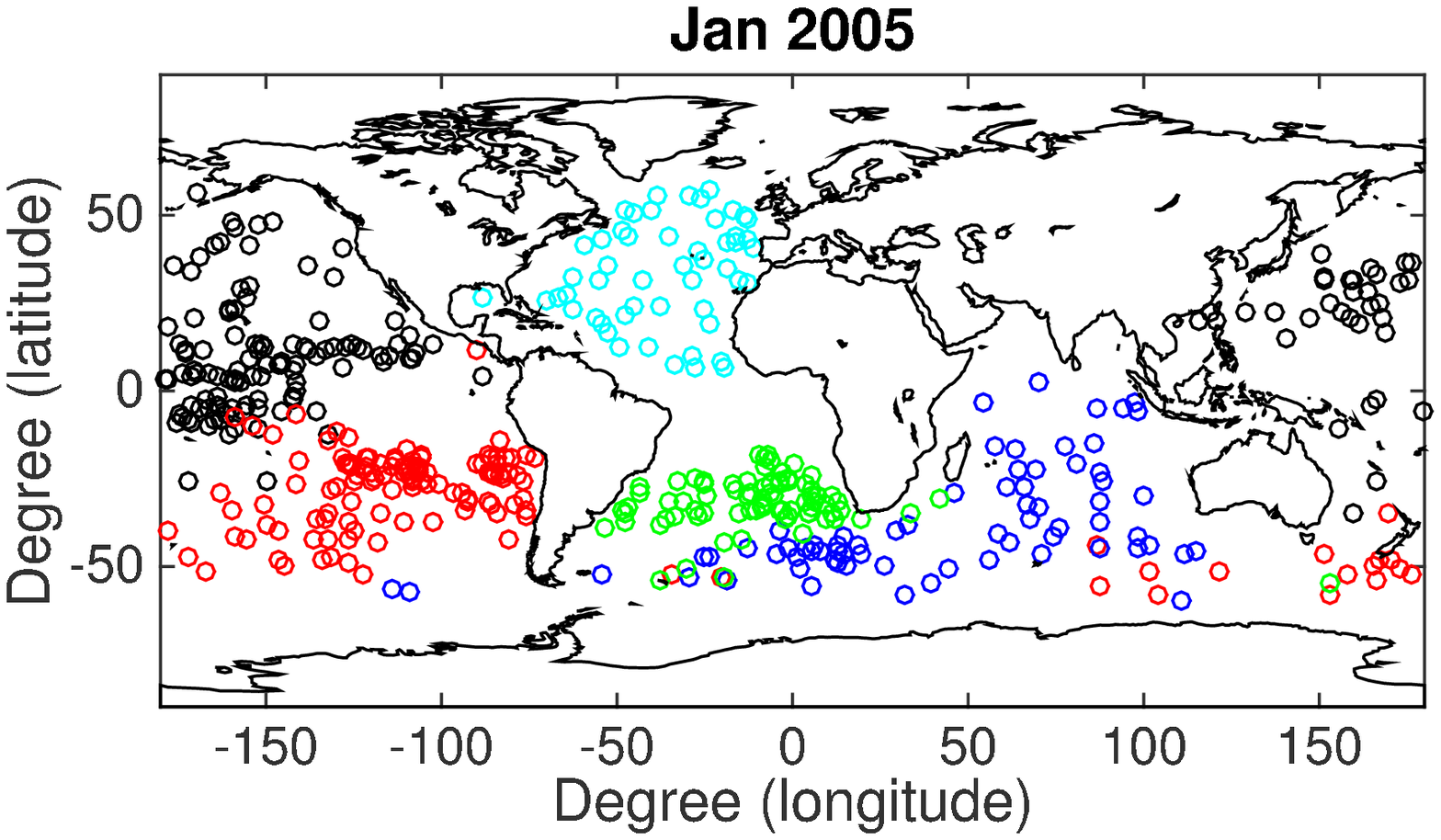} \\
{\scriptsize (a)}\\[2mm]
\includegraphics[width=0.98\columnwidth, clip=true]{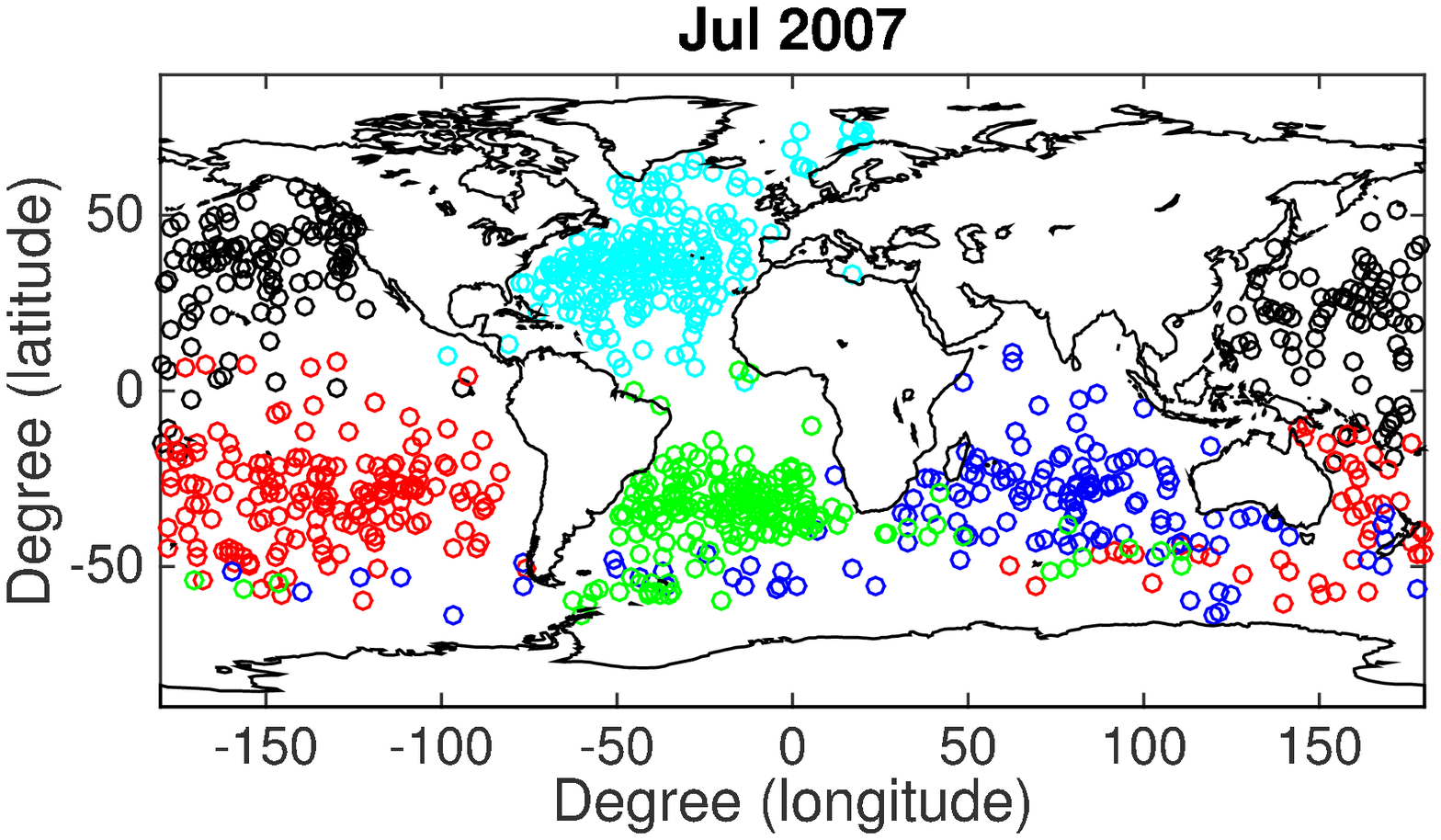}\\
{\scriptsize (b)}\\[2mm]
\includegraphics[width=0.98\columnwidth, clip=true]{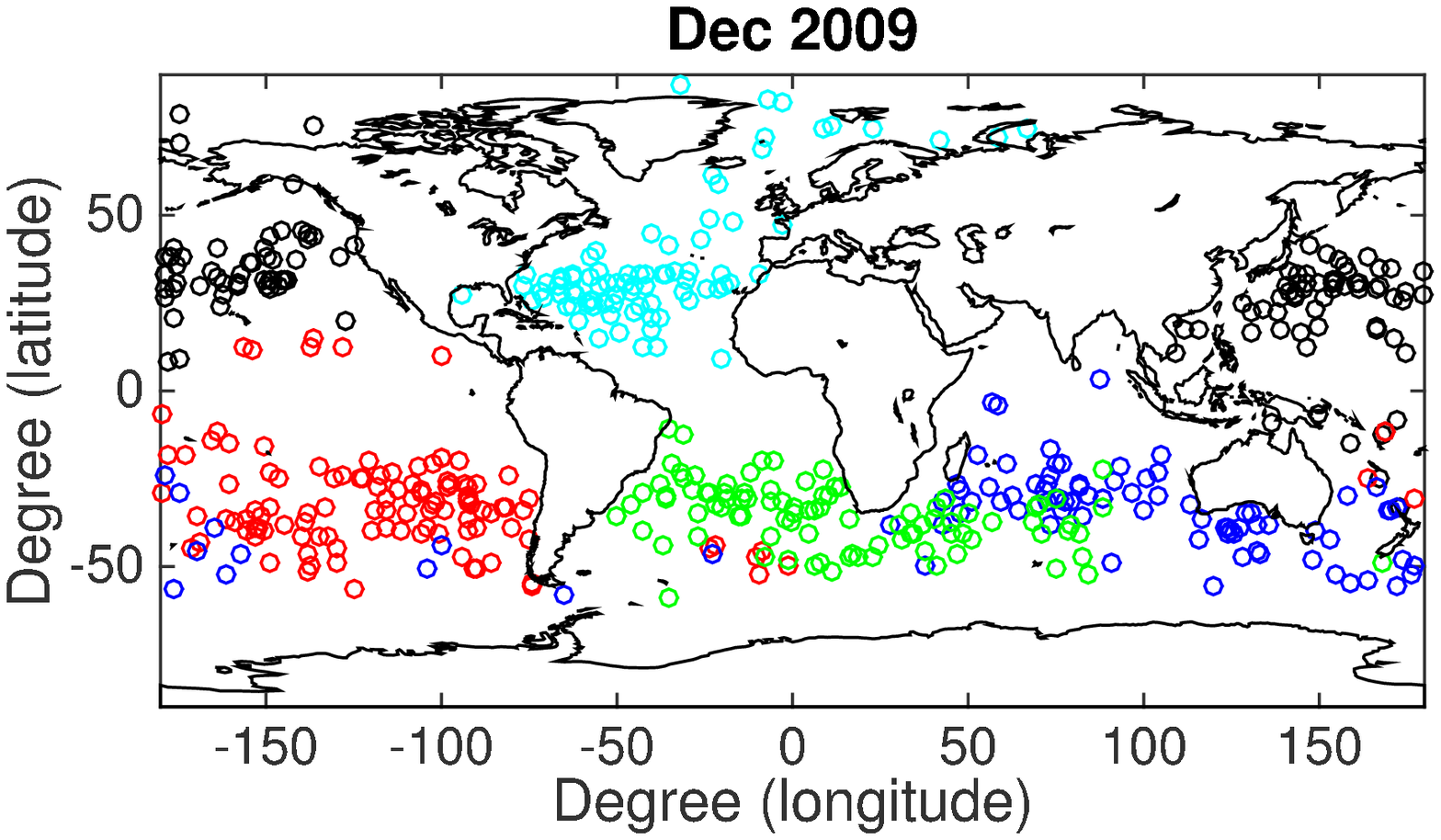}\\
{\scriptsize (c)}\\
\end{tabular}
\caption{5-clustering of drifter data ($m=1.5$). \\(a) drifter positions January 2005, (b) drifter positions July 2007, (c) drifter positions December 2009. Animation of the cluster motion 2005--2009 in online version (Multimedia view).}\label{fig:drifter5}
\end{figure}

Figure \ref{fig:drifter5} shows the results of Algorithm 2 with $K=5$ at January 2005 (a), July 2007 (b), and December 2009 (c). An animation is available online (Figure \ref{fig:drifter5} (Multimedia view)).
We choose $K=5$ in order to attempt to delineate the five major oceans:  the North and South Atlantic Oceans, the North and South Pacific Oceans, and the Indian Ocean.
Broadly, we see that the clustering does find the appropriate equatorial separations of the Atlantic and Pacific Oceans, and also separates the Indian Ocean.

\begin{figure}[h!]
\includegraphics[width=0.98\columnwidth, clip=true]{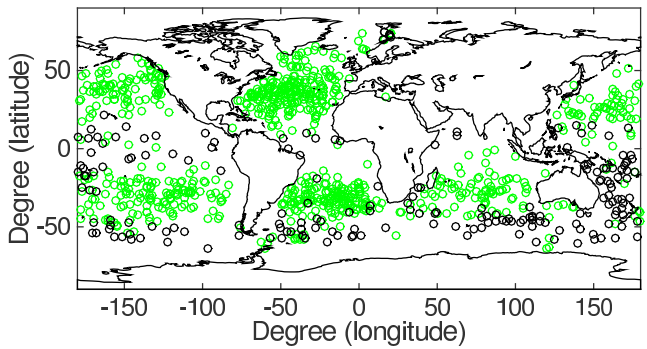} \\
\caption{%Uncertainty of cluster membership ($m=1.5$, $K=5$) using entropy (\ref{eq:entropy}).
Difficult to classify drifters (positions as of July 2007) are marked black when their relative entropy is $h>0.1$. Animation of the cluster motion 2005--2009 in online version (Multimedia view).}\label{fig:drifter5_entrop}
\end{figure}

Some of these separations are highlighted by investigating the certainty of membership of the individual drifters based on an entropy calculation (\ref{eq:entropy}).
In Figure \ref{fig:drifter5_entrop} those drifters (positions as of July 2007) are marked black when their relative entropy is $>0.1$, corresponding to a maximum membership value of less than $\approx 0.96$.
Figure \ref{fig:drifter5_entrop} and in particular the time evolution of the drifters (Figure \ref{fig:drifter5_entrop} (Multimedia view)) shows that the uncertain regions correspond to the major ocean barriers in the Atlantic and Pacific, and the Southern Ocean.

Our results in Figure \ref{fig:drifter5} are strikingly similar to those shown in Figure 6 in Ref.~\onlinecite{FSvS14}, which have been derived using transfer operator methods and the (wind-forced) OFES model.
For example in Figure \ref{fig:drifter5} (a), when comparing with Figure 6 in Ref.~\onlinecite{FSvS14} we see: the separation of the Pacific Ocean becoming more southerly as one proceeds westwards toward Australia;  the Indian Ocean spilling westwards at its southerly boundary; and the South Atlantic forcing its way around the east coast of southern Africa.
As described in Ref.~\onlinecite{FSvS14}, the Ekman dynamics of the ocean surface circulation guarantees that each of the five major oceans contains an attracting region associated to the great oceanic gyres and their corresponding garbage patches.
The separations seen in Figure \ref{fig:drifter5} (a) and Figure 6 in Ref.~\onlinecite{FSvS14} and the features described above are associated with the basins of attraction of these five attracting regions.

We remark that while the results in Ref.~\onlinecite{FSvS14} are of a higher spatial resolution than those obtained here, the experiments in Ref.~\onlinecite{FSvS14} used just over $10^6$ trajectories, recorded every eight weeks for a period of 48 weeks (a total of $6.14\times 10^6$ data points), while here we have $869\times 60=5.2\times 10^4$ data points, comprised of 2267 incomplete trajectories.
We also remark that while we have drawn comparisons between Ref.~\onlinecite{FSvS14} and the present study, the former computed ocean boundaries as basins of attraction based on a repeating 48-week ocean circulation, while our present study seeks to compute estimates of coherent sets based on five years of non-repeating drifter data.

\section{Discussion}
We have introduced a ``rough-and-ready'' general cluster-based approach for analysing coherent structures in time-dependent dynamical systems.
Our method assigns individual trajectories membership in regions that retain a compact extent over a specified finite time duration.
Our approach has several advantages.

First, the ability to work directly with a small number of trajectories, including the situations where the trajectories do not span the entire time duration of interest and where observations may be missing from within trajectories.
Second, initial implementation is rapid (using e.g.\ the built-in MATLAB function \verb"fcm" to perform the fuzzy clustering for the case of complete data), and the runtimes are fast (on the order of fractions of seconds for the one-dimensional maps in Section \ref{sec:1d} to less than 10 seconds to cluster a dataset of $32768$ trajectories in $202$ dimensions, as in the case of the double gyre flow with flow duration $\tau=10$ in Section \ref{sec:dg}).
Third, our method considers entire trajectories (not just the endpoints) and automatically outputs clusters at every time instant in the trajectory data;  thus a frame-by-frame description of the temporal evolution of the clusters is immediately obtained.
Fourth, the use of fuzzy clustering provides feedback in the form of membership likelihoods and entropy, which provide the user with an estimate of confidence with which a trajectory has been assigned to a particular compact region.
Finally, the soft clustering approach is relatively insensitive to noise in the data.
We note that the same methodology can be used to estimate coherent regions for SDEs, by simply generating stochastic trajectories and applying Algorithm 1.
%Future improvements could include a refinement of the clustering method to include connectivity or density elements in order to move away from favouring spherical clusters of similar size.
%One possibility is DBSCAN\cite{Ester_et_al1996}, which does not require one to fix the number of clusters beforehand and is not restricted to spherical clusters.

\begin{acknowledgments}
%GF and KPG thank Michael Allshouse for communicating the double gyre example EQUATION NUMBER HERE.
The research of GF is supported by an ARC Future Fellowship (FT120100025).  GF also thanks the University of Canterbury's Erskine Fellowship scheme for partial financial support and the Department of Mathematics and Statistics at the University of Canterbury for providing excellent working conditions during part of the time this research was undertaken. KPG acknowledges support from an ARC Discovery Project (DP110100068). She also thanks the School of Mathematics and Statistics at the University of New South Wales for hospitality as well as the University of Canterbury for hospitality and financial support. \end{acknowledgments}

\vspace*{5mm}
\noindent\textbf{REFERENCES}

\bibliographystyle{abbrv}
\bibliography{cluster}

\end{document}